\documentclass[10pt,a4paper,reqno]{amsart}
\numberwithin{equation}{section}
\usepackage{amsmath,amssymb,amsthm}
\usepackage{xcolor,textcomp}
\usepackage{graphicx, epstopdf}
\usepackage{braket,amsfonts}
\usepackage{dsfont}
\usepackage{mathtools}
\usepackage{mwe}
\usepackage[T1]{fontenc}
\usepackage{xfrac}

\usepackage[top=3cm,bottom=3cm,left=3cm,right=3cm]{geometry}

\usepackage{hyperref}
\hypersetup{colorlinks=true,linkcolor=red,citecolor=blue,filecolor=magenta,urlcolor=cyan} 
\usepackage[capitalise]{cleveref}
\usepackage{comment}
\usepackage{tikz}
\usepackage{pgfplots}
\usepackage{stmaryrd}
\usepackage[mathscr]{euscript}
\usepackage[sort,nocompress]{cite}

\newtheorem{theorem}{Theorem}[section]
\newtheorem{lemma}[theorem]{Lemma}
\newtheorem{remark}[theorem]{Remark}
\newtheorem{coro}[theorem]{Corollary}
\newtheorem{defi}[theorem]{Definition}
\newtheorem{prop}[theorem]{Proposition}

\usepackage{amsopn}
\DeclareMathOperator*{\esssup}{\textnormal{ess\,sup}}

\DeclareMathOperator{\weak}{\textnormal{weak}}

\DeclareMathOperator{\tr}{tr}

\DeclareMathOperator{\inn}{\textnormal{in}}
\DeclareMathOperator{\out}{\textnormal{out}}
\DeclareMathOperator{\sym}{\textnormal{sym}}

\newcommand\C{\mathcal{C}}

\newcommand\R{\mathbb{R}}

\newcommand\Ss{\mathbb{S}}

\newcommand\eps{\epsilon}
\newcommand\veps{\varepsilon}

\newcommand{\boldvphi}{\boldsymbol{\vphi}}  

\newcommand{\disp}{\displaystyle}

\newcommand{\vc}[1]{{\bf #1}}

\newcommand{\uvect}{\mathbf{u}}

\newcommand{\Div}{\textnormal{div}}
\newcommand{\D}{\mathbb{D}}

\newcommand{\F}{\mathcal{F}}

\newcommand{\vphi}{\varphi}
\renewcommand{\rho}{\varrho}

\def\dx{\,\textnormal{d}x}
\def\dt{\textnormal{d}t}
\def\d{\,\textnormal{d}}

\makeatletter
\@namedef{subjclassname@2020}{%
	\textup{2020} Mathematics Subject Classification}
\makeatother

\title[Compressible active nematic liquid crystals]{Dissipative solutions to a Beris-Edwards type model for compressible active nematic liquid crystals}

\author[K. Bhandari, A. Majumdar, {\v{S}}. Ne{\v{c}}asov{\'{a}}]{Kuntal Bhandari$\, ^{1}$ 
 \and Apala Majumdar$\, ^{2}$	
\and {\v{S}}{\'{a}}rka 
Ne{\v{c}}asov{\'{a}}$\, ^{3}$ }

\thanks{$^{1}$Institute of Mathematics of the  Czech Academy of Sciences, \v{Z}itn\'a 25, 11567 Praha 1, Czech Republic; \texttt{bhandari@math.cas.cz}.}

\thanks{$^2$Department of Mathematics, University of Manchester, Manchester M13 9PL, United Kingdom; 
\texttt{apala.majumdar@manchester.ac.uk}.}

\thanks{$^{3}$Institute of Mathematics of the  Czech Academy of Sciences, \v{Z}itn\'a 25, 11567 Praha 1, Czech Republic;  \texttt{matus@math.cas.cz}. }



\date{\today}

\begin{document}

\begin{abstract} 

We study the  hydrodynamics of compressible active nematic liquid crystals  in a three-dimensional and bounded domain, with a  nonlinear viscosity tensor and nonhomogeneous boundary data, in a Landau-de Gennes framework. We prove the existence of  {\em dissipative solutions}  within a Beris-Edwards type model for active nematodynamics,  which are  weak solutions satisfying the underlying equations modulo a defect measure. The proof follows from a three level approximation scheme -- the Galerkin approximation, the classical parabolic regularization of the continuity equation, and the convex regularization  of the potential generating the viscous stress. New techniques are required to deal with non-Newtonian stress tensor, larger classes of admissible pressure potentials and nonhomogeneous boundary conditions. 

    
\end{abstract}

\maketitle

\tableofcontents

\section{Introduction}

Nematic liquid crystals (NLCs) are classical examples of partially ordered materials that combine fluidity with the ordering of crystalline solids \cite{DeGennes}. NLCs are anisotropic complex fluids with constituent rod-like or anisotropic molecules. The rod-like molecules tend to align along locally preferred directions or distinguished material directions referred to as \emph{nematic directors}, and NLCs consequently have long-range orientational order and direction-dependent physical, optical and rheological properties. NLCs have long been celebrated as the working material of choice for the multi-billion dollar liquid crystal display (LCD) industry and various electro-optic and photonic devices. Active nematic liquid crystals are a relatively new and rapidly growing scientific field, relevant for biological systems, biomimetic materials and new classes of animate and autonomous materials \cite{majumdarphilosophicaltransactions2014}. Active nematics typically correspond to systems composed of elongated particles that exhibit spontaneous coordinated and collective motion, i.e., they are fundamentally non-equilibrium systems that are constantly driven out of equilibrium by internal energy sources \cite{giomi-mahadevan-chakraborty-hagan-2012}. Some examples of active systems are swarms of micro-organisms or bacterial colonies \cite{darnton2010dynamics}, systems of microtubules \cite{sanchez2012spontaneous} and microswimmers \cite{wensink2012meso}. Active nematics have exceptional properties such as giant density fluctuations \cite{mishra2010fluctuations,narayan2007long}, spontaneous laminar flows \cite{giomi2008complex,marenduzzo2007steady}, low Reynolds number turbulence \cite{wensink2012meso,giomi-mahadevan-chakraborty-hagan-2012}, exotic spatial and temporal textures with fascinating defect patterns \cite{ginelli2010large,saintillan2008instabilities} to name a few. As such, there is a burgeoning interest to rigorously analyse theoretical frameworks for active nematics in various settings and we present the first rigorous analysis of \emph{dissipative solutions} in this context, a relatively novel class of solutions that are wider in scope than conventional weak solutions and much needed for comprehensive analyses of numerical methods for complex flows.

There are different competing continuum or macroscopic theories for NLCs: the simplest Oseen-Frank theory \cite{Frank1958liquid,Oseen1933theory} restricted to uniaxial NLCs with a single distinguished nematic director with constant scalar order parameter; the Ericksen theory \cite{Ericksen1961conservation} restricted to uniaxial nematics with variable scalar order parameter and the most general Nobel-Prize winning Landau-de Gennes (LdG) theory \cite{DeGennes} that can account for NLC phases with multiple directors and order parameters, to account for generic NLC phenomena and complex NLC defects \cite{liustaticdynamictheory}. The order parameters measure the degree of orientational ordering about the nematic directors and the defect set is associated with the nodal set of the uniaxial order parameter in the Ericksen theory. This is a restrictive definition of nematic defects and misses important phenomena such as biaxiality (a primary and secondary nematic director) and \emph{escape into third dimension} near defect sets. We work in the powerful and most general LdG framework for NLCs, which can be generalised to active nematics. In the LdG framework, the NLC state is described by a macroscopic order parameter, the ${Q}$-tensor order parameter, which can be defined in terms of experimental quantities such as the dielectric anisotropy or optical birefringence data \cite{DeGennes,maityproceedingsroyalsociety2025}. Mathematically, the ${Q}$-tensor is a symmetric, traceless $3\times 3$ matrix whose eigenvectors model the nematic directors and the corresponding eigenvalues define the order parameters about the nematic directors. The ${Q}$-tensors belong to the space
\begin{align*}
S_0^3 = \left\{Q \in \mathbb R^{3\times 3} \mid Q= Q^\top , \ \tr(Q)=0   \right\}. 
\end{align*} The NLC is said to be in the isotropic phase if ${Q}=0$, uniaxial phase if ${Q}$ has two degenerate non-zero eigenvalues and in the biaxial phase if ${Q}$ has three distinct eigenvalues. In particular, a uniaxial ${Q}$-tensor can be written as 
\[
{Q} = s\left(\mathbf{j}\otimes\mathbf{j} - \frac{1}{3}\mathbb{I}_3 \right)
\]
where $\mathbf{j}$ is the uniaxial nematic director or the eigenvector with the non-degenerate eigenvalue and $s$ is the associated scalar order parameter.
The LdG theory is a variational theory and the LdG free energy is a nonlinear and non-convex functional of ${Q}$ and $\nabla_x {Q}$ as given below:
\begin{equation}\label{eq:LdGenergy}
\F = \int_\Omega \left(  \frac{K}{2} |\nabla_x Q|^2  +  \frac{k}{4} (c-c_*) \tr(Q^2)  - \frac{b}{3} \tr(Q^3) + \frac{c_*}{4} \tr^2(Q^2)  \right).
\end{equation}
Here we use the one-constant elastic energy density, $K$ is a positive elastic constant, $c_*$ is the critical concentration for the isotropic-nematic transition, and $k>0$ and $b\in \R$ are material-dependent constants \cite{majumdar2010,C-Majumdar-W-Z-Comp}. Without loss of generality, we take $K=k=1$. The equilibrium configurations are mathematically modelled by minimisers of the LdG theory subject to the imposed boundary conditions.

Active nematics are intrinsically non-equilibrium fluids, often thought of as a suspension of elongated active units in an ambient nematic fluid. We work with an active version of the compressible Beris-Edwards model for nematodynamics \cite{berisedwards}. The Beris-Edwards theory models the time evolution of the LdG ${Q}$-tensor and the fluid velocity, with non-linear couplings between the two evolution equations. The evolution of the velocity is dictated by a Navier-Stokes type equation with nonlinear stresses that originate from the LdG ${Q}$-tensor, which makes the problem significantly harder than standard fluid mechanics problems. The evolution of the LdG ${Q}$-tensor also contains coupling terms between the LdG order parameter and the velocity gradient, leading to new technical difficulties. The active version of the Beris-Edwards model has four model parameters: the concentration $c$ of the active particles, the density $\rho$ of the homogenized fluid, the fluid velocity $\uvect$, and the LdG ${Q}$-tensor order parameter. Following \cite{C-Majumdar-W-Z-Comp}, we employ the following system of partial differential equations to model compressible active nematic fluids in three dimensions.

Consider a bounded Lipschitz (spatial) domain $\Omega \subset \mathbb R^3$ and time interval $(0,T)$. Then the time evolution of $\left(\rho, \uvect, c, {Q}\right)$ is given by:
\begin{align}
& \partial_t \rho  + \Div_x (\rho \uvect) = 0, \label{dens-eq} \\
& \partial_t (\rho \uvect) + \Div_x (\rho \uvect \otimes \uvect) + \nabla_x p(\rho) = \Div_x \Ss + \Div_x(\tau + \sigma)  , \label{momen-eq} \\
& \partial_t c + (\uvect \cdot \nabla_x) c = D_0 \Delta_x c    , \label{concen-eq}  \\
& \partial_t Q + (\uvect \cdot \nabla_x) Q + Q \Lambda - \Lambda Q = \Gamma H[Q, c] , \label{nematic-tensor-eq} 
\end{align} 
where $D_0>0$ is the diffusion constant and  $\Lambda : = \frac{1}{2}\left(\nabla_x \uvect - \nabla_x^\top \uvect \right)$ is the antisymmetric part of the strain tensor. Furthermore, the  pressure $p=p(\rho)$  depends on the fluid density and its precise properties are given later in Section \ref{Section-Pressure}.
Motivated from  \cite[Section 2.1.2]{Dissipative-Feireisl}, we assume that the viscous stress tensor $\Ss$ is related to the symmetric velocity gradient 
$$\D_x\uvect = \frac{1}{2}\left(\nabla_x \uvect + \nabla_x^\top \uvect\right)$$
through a {\em general implicit rheological law}
\begin{align}\label{stress-tensor}
\Ss : \D_x \uvect = F(\D_x\uvect) + F^*(\Ss) , 
\end{align}
for a suitable potential  $F$ (with its conjugate $F^*$) such that 
\begin{align}\label{F-convex}
F : \R^{3\times 3}_{\sym} \to [0, \infty) \text{ is a  convex lower semi-continuous (l.s.c.) function with } F(0)=0, 
\end{align}
where $\R^{3\times 3}_{\sym}$ is the space of three-dimensional real symmetric tensors. 
Since $F$ is proper convex l.s.c. function, the conjugate  $F^* : \R^{3\times 3}_{\sym}\to [0,\infty)$ is also convex, l.s.c. as well as superlinear, cf. \cite[Proposition 2.1]{Basaric}. Moreover, due to the proper convexity and lower semicontinuity of the function $F$,  the relation  \eqref{stress-tensor} can be interpreted in view of Fenchel--Young inequality, given by 
\begin{align}\label{relation-stress-tensor}
\Ss \in \partial F(\D_x \uvect) \Leftrightarrow  \D_x \uvect \in \partial F^*(\Ss),
\end{align}
where $\partial$ denotes the subdifferential/ subgradient  of a convex function.

In this work, we particularly make the following assumption on $F$, namely   
\begin{align}\label{Coercivity-F} 
F(\D) \geq \mu_0 \left| \D - \frac{1}{3} \tr [\D] \mathbb I_3 \right|^{\sfrac{4}{3}} \quad \text{for all } \, \D \in \mathbb R^{3\times 3}_{\sym}  ,
\end{align}
for some $\mu_0>0$, and the choice of the power $\frac{4}{3}$ is specific to the active stresses in the evolution equation for $\uvect$. Here and in the sequel --  for any matrix $A : = (a_{ij})_{1\leq i,j\leq 3} \in \R^{3\times 3}$, we denote by $|A|$, the Frobenius  norm or Hilbert-Schmidt norm of $A$ over the field $\mathbb R$, that is 
\begin{align*}
|A| := \bigg(\sum_{1\leq i,j \leq 3} |a_{ij}|^2 \bigg)^{\sfrac{1}{2}} = \left[\tr(A^\top A)\right]^{\sfrac{1}{2}} . 
\end{align*} 
Note that, for the standard Newton's rheological law, we have
\begin{align*}
\Ss = \mu \D_x \uvect + \lambda \Div_x \uvect \mathbb I_3 ,
\end{align*}
and accordingly, the Newtonian viscous stress tensor is then given by
\begin{align*}
F(\D_x \uvect) = \frac{\mu}{2} |\D_x \uvect|^2 + \frac{\lambda}{2} |\Div_x \uvect \mathbb I_3 |^2 , \ \ \mu >0 , \ \lambda + \frac{2}{3} \mu \geq 0 .
\end{align*}
The resulting problem is the  compressible system for nematic liquid crystals with standard Newton's rheological law;  see,  for instance \cite{C-Majumdar-W-Z-Comp}. In particular, we consider a wider class of non-Newtonian viscous stress tensors in this paper.

The stress tensor $\sigma$ has the following form in \eqref{momen-eq}:
\begin{align*}
\sigma = \sigma_{r} + \sigma_a ,
\end{align*}
with 
\begin{align*}
\sigma_r = Q H[Q,c] - H[Q,c]Q,  
\end{align*}
and 
\begin{align}\label{sigma-a} 
{\sigma_{a} =\sigma_* c^2 Q} ,
\end{align}
where $\sigma_r$ is the stress due to the nematic elasticity and $\sigma_a$ is the {\bf active stress} which models contractile ($\sigma_*>0$) or extensile ($\sigma_*<0$) stresses exerted by the active particles along the director field.

In equations \eqref{momen-eq} and \eqref{nematic-tensor-eq},  $\Gamma^{-1}>0$ is the rotational viscosity and $H[Q,c]$ is given  by 
\begin{align*}
H[Q, c] : =  K \Delta_x Q - \frac{k}{2} (c-c_*) Q + b \left( Q^2 - \frac{\tr (Q^2)}{3} \mathbb I_3 \right) - c_* Q \tr(Q^2) ,  
\end{align*}
which describes the relaxation  dynamics in terms of the gradient flow of the LdG free energy, i.e., $H_{\alpha \beta} = \frac{-\delta \F}{\delta Q_{\alpha \beta }}$ where $F$ has been defined in \eqref{eq:LdGenergy}.

The symmetric additional stress tensor $\tau$ is denoted by 
\begin{align*}
\tau=  G(Q) \mathbb I_3 - \nabla_x Q \odot \nabla_x Q ,
\end{align*}
where
\begin{align*}
G(Q) = \frac{1}{2} |\nabla_x Q|^2 + \frac{1}{2} \tr(Q^2) + \frac{c_*}{4} \tr^2(Q^2).
\end{align*}

Given these definitions, we rewrite the system \eqref{dens-eq}--\eqref{nematic-tensor-eq}  explicitly as follows:
\begin{align}
& \partial_t \rho  + \Div_x (\rho \uvect) = 0, \label{dens-eq-1} \\
& \partial_t (\rho \uvect) + \Div_x (\rho \uvect \otimes \uvect) + \nabla_x p(\rho) = \Div_x \Ss + 
\Div_x \left(  G(Q) \mathds{I}_3 - \nabla_x Q \odot \nabla_x Q  \right) \notag \\
& \qquad \qquad \qquad \qquad \qquad \qquad \qquad \qquad  + \Div_x \left(Q \Delta_x Q - \Delta_x Q Q     \right) + {\sigma_* \Div_x (c^2 Q)} , \label{momen-eq-1} \\
& \partial_t c + (\uvect \cdot \nabla_x) c = D_0 \Delta_x c, \label{concen-eq-1} \\
& \partial_t Q + (\uvect \cdot \nabla_x) Q + Q \Lambda - \Lambda Q = \Gamma H[Q, c] , \label{nematic-tensor-eq-1} 
\end{align} 
where 
\begin{align}\label{formula-H}
H[Q,c] = \Delta_x Q - \frac{1}{2} (c-c_*) Q + b \left( Q^2 - \frac{\tr (Q^2)}{3} \mathbb I_3 \right) - c_* Q \tr(Q^2)  .
\end{align}
 The evolution equations are supplemented with the following initial and boundary conditions.

\vspace{.2cm}
\noindent
\underline{\em Initial conditions.} We consider 
\begin{align}\label{initial-concen}
&c(0)=c_0 \in L^\infty(\Omega), \ 0< \underline{c} \leq c_0 \leq \overline{c} <\infty \\
&\rho (0) = \rho_0 \in L^\gamma(\Omega), \ \gamma > 1,  \ \rho_0\geq 0 , \label{initial-dens} \\
&(\rho \uvect)(0) = (\rho \uvect)_0 , \  (\rho \uvect)_0 = 0 \text{ if } \rho_0 = 0 , 
\  \frac{|(\rho \uvect)_0|^2}{\rho_0} \in L^1(\Omega), \label{initial-velocity} \\
& Q(0)=Q_0 \in H^1(\Omega), \ Q_0 \in  S^3_0 \, \text{ a.e. in } \Omega ,\label{initial-Q} 
\end{align}
where $\underline{c}$, $\overline{c}$ are constants. 

\vspace{.2cm}
\noindent 
\underline{\em Boundary conditions.}
For the density and fluid velocity, we consider general inflow-outflow boundary conditions: 
\begin{align}\label{inflow-outflow} 
\uvect |_{\partial \Omega} = \uvect_B \ \text{ and } \ \rho|_{\Gamma_{\inn}} = \rho_B ,
\end{align}
where 
\begin{align*}
\ \Gamma_{\inn} = \left\{ x \in \partial \Omega \mid  \uvect_B(x) \cdot \mathbf n(x) < 0 \right\},  
\end{align*}
 $\mathbf n$ denotes the unit outward normal to  $\partial \Omega$. Accordingly, the boundary $\partial \Omega$ is decomposed as 
\begin{align*}
\partial \Omega = \Gamma_{\inn} \cup \Gamma_{\out} ,
\end{align*} 
where the inflow part $\Gamma_{\inn}$ is given above and the outflow part is 
\begin{align*}
\Gamma_{\out} = \left\{ x \in \partial \Omega \mid  \uvect_B(x) \cdot \mathbf n(x) \geq  0 \right\}.  
\end{align*}
 For the sake of simplicity,  one may consider that $\uvect_B\in \C^1_c(\mathbb R^3; \mathbb R^3)$ (that is, $\uvect_B$ can be extended smoothly throughout the whole space) and $\rho_B \in \C^1(\partial \Omega)$ with $\rho_B \geq \overline{\rho}>0$.

The usual boundary conditions for the LdG ${Q}$-tensor order parameter are the homogeneous Neumann or   (non)homogeneous Dirichlet conditions, see for instance \cite{Paicu-Zarnescu-ARMA,C-Majumdar-W-Z-Comp,Abels-Dolzman-Liu,Rodriguez-Bellido-Gongalez,Wang-Xu-Yu}, and references therein. In the present paper, motivated from \cite{Rodriguez-Bellido-Gongalez},  we impose the following nonhomogeneous Dirichlet boundary condition,  
\begin{align}\label{boundary-Q}
Q |_{\partial \Omega} : = Q_B \in H^{\sfrac{3}{2}}(\partial \Omega).
\end{align}
Finally, we  impose a homogeneous Neumann  boundary condition for the concentration of active particles, i.e., 
 \begin{align}\label{boundary-c}
\nabla_x c \cdot \mathbf n = 0 \ \text{ on } (0,T) \times \partial \Omega  .
 \end{align} 
 Note that the above boundary condition \eqref{boundary-c} is also considered in \cite{C-Majumdar-W-Z-Comp}. Further, the  boundary data $\uvect_B$, $\rho_B$ and $Q_B$ are independent of time and only depend on the spatial coordinates of 
 $\partial \Omega$. 

 We briefly review existing work on nematodynamics, although our references are by no means exhaustive. The system of partial differential equations \eqref{dens-eq-1}--\eqref{nematic-tensor-eq-1} reduces to passive nematodynamics if the active stress is zero and to the incompressible case if the density is constant. In \cite{Paicu-Zarnescu-SIMA}, the authors prove the existence of weak solutions for the incompressible passive Beris-Edwards model for nematodynamics, in two and three dimensions. In \cite{Wang-Xu-Yu}, the authors study a compressible model for coupled nematodynamics and prove the existence of global weak solutions for the same. In \cite{C-Majumdar-W-Z-InComp}, the authors prove the existence of a global weak solution for the incompressible version of the active Beris-Edwards model, including active stresses, in two and three dimensions. There are parallel results for the 
 Ericksen-Leslie theory for incompressible flows of uniaxial nematics, see \cite{Hieber-Pruss-Incomp} (and references therein), and for work on compressible uniaxial nematic flows, see\cite{Hieber-Pruss-Comp}. We also mention the work of \cite{Geng-Roy-Arghir,Binz-Roy-Hieber-Brandt}, wherein the authors analyse the liquid crystal--rigid body interaction models for incompressible flows.  
 Our work builds on the work in \cite{C-Majumdar-W-Z-Comp} wherein the authors study the  compressible (Newtonian) version of the active Beris-Edwards model for nematodynamics, i.e., they prove the existence of a {\em weak solution} for the system \eqref{dens-eq-1}--\eqref{nematic-tensor-eq-1} on bounded, three-dimensional domains subject to prescribed initial and boundary conditions. The system \eqref{dens-eq-1}--\eqref{nematic-tensor-eq-1} is not the full model for compressible active nematodynamics -- it assumes isotropic diffusion for the concentration equation; the flow-aligning parameter $\lambda$ is set to zero in the full set of equations for compressible active nematodynamics proposed in \cite{giomi-mahadevan-chakraborty-hagan-2012} which simply means that we are in the flow-tumbling regime, and we neglect one of the terms in the passive nematic stress. However, this system captures essential salient features of active nematodynamics, under suitable hypotheses. In \cite{C-Majumdar-W-Z-Comp}, the authors use a Faedo-Galerkin three-step approximation argument to prove the existence of weak solutions for the system \eqref{dens-eq-1}--\eqref{nematic-tensor-eq-1}: (i) an artificial pressure approximation, (ii) artificial viscosity approximation and (iii) passage from finite-dimensional to infinite-dimensional spaces. There are multiple technical challenges stemming from the nonlinear coupling stresses between ${Q}$ and $\mathbf{u}$, and notably the active stress tensor, $\sigma_a$. The authors introduce new technical tools to circumvent these challenges and exploit the maximum principle for the concentration equation, the nature of the stress tensor $\tau$ and the symmetry and tracelessness properties of ${Q}$ to obtain judicious cancellations and well-behaved inequalities for the mathematical analysis of the same.  

 Our work differs from the work in \cite{C-Majumdar-W-Z-Comp} in three important ways. We prove the existence of \emph{dissipative solutions} for the compressible model of active nematodynamics in \eqref{dens-eq-1}--\eqref{nematic-tensor-eq-1} with inhomogeneous boundary data for the LdG ${Q}$-tensor and the flow velocity, along with inflow boundary conditions for the density. Dissipative solutions have been introduced  in \cite{Feireisl-Abbatielo} for general models of incompressible fluids, as a larger class of general solutions than the standard weak solutions. Dissipative solutions for the  compressible Euler equations  have been studied  in \cite{Feireisl-Breit}. For dissipative solutions, the viscous stress tensor $\Ss$ is replaced by a {\bf Effective viscous stress} $\Ss_{\textnormal{eff}}: = \Ss - \Re$, where $\Re$ is called the {\bf Reynolds stress}. Dissipative solutions reduce to classical weak solutions if the Reynolds stress $\Re =0$ and thus one can think that $\Re$  as an error term in the momentum equation.  Accordingly, this generates another error term  $\mathfrak{E}$ in the energy inequality, referred to as the {\bf Energy dissipation effect}. Dissipative solutions have attracted substantial interest in the applied analysis community, for example, 
 dissipative solutions for 
compressible viscous fluids are  rigorously  studied  in \cite{Dissipative-Feireisl} with a pressure law  of type $p(\rho)=a\rho^\gamma$ with $\gamma>1$  (see \eqref{S1} of this paper).  A similar  result with a linear pressure law is analysed in \cite{Basaric} (the isentropic case with the adiabatic constant $\gamma=1$). In \cite{Sarka-Jan-Justyna}, the authors address the existence of   dissipative solutions for compressible Navier-Stokes systems with Coulomb friction law boundary condition. 

We successfully  extend the machinery in \cite{Feireisl-Abbatielo} to active nematodynamics for the first time. There are some immediate advantages of this approach. One cannot prove the existence of standard weak solutions for compressible fluids with the pressure law $p(\rho)=a\rho^\gamma$ and adiabatic constant  $\gamma>1$; we need the condition $\gamma>\frac{3}{2}$ for that purpose. However, one can allow for $\gamma>1$ in the pressure law whilst studying {\bf dissipative solutions} for compressible fluids. The work in \cite{C-Majumdar-W-Z-Comp} considers a Newtonian viscous stress but we allow for more general rheological laws and viscous stress tensors, that could be non-Newtonian in principle. In fact, one can relax the condition \eqref{Coercivity-F}  if we eliminate the active stress and focus on compressible passive nematodynamics (see Section \ref{Section-Conclusion} for further details).  We also allow more general boundary data in our formulation, namely the so-called inflow boundary condition for the density and nonhomogeneous boundary data for the velocity $\uvect$ and for the $Q$-tensor parameter. The inflow/outflow boundary conditions are important in many real-world applications. In fact this is a natural setting for blood flow in arteries,  flows in pipelines, wind tunnels and turbines; see for example \cite{Adelia-1,Blood-flow-1,Mohammadi2009boundary} and references therein. For compressible flows in barotropic case, the existence of global-in-time weak solutions with general inflow/outflow boundary conditions has been analyzed in \cite{ChJiNo,Choe-Nov-Yang}. We also refer to \cite{Novo,Girinon} where such boundary conditions are considered but with several restrictions on the boundary and boundary data.  With regards to strong solutions for the compressible Navier-Stokes system with the boundary conditions \eqref{inflow-outflow}, we quote \cite{Mucha-Zazackowski} and some references therein. Finally, in \cite{Jan-Valasek}, the authors carry out a numerical study of an incompressible fluid-structure interaction problem with inflow type boundary condition. Hence, these boundary conditions significantly enhance the potential practical implications of our study. 
  

\subsection*{Paper organization}
The rest of the paper is organized as follows. In Section \ref{Sec-weak-formulation}, we define the concept of dissipative solutions and state our main existence result.
Section \ref{Section-Apriori} is devoted to a-priori estimates which are crucial for the derivation of the energy inequality in our model. The remaining sections are devoted to prove the existence of dissipative solutions for the compressible system \eqref{dens-eq-1}--\eqref{nematic-tensor-eq-1}, and the proof consists of several approximation schemes and limit passages. More precisely, 
in Section \ref{Section-Existence}, we prove the existence of approximate  solutions associated with  our system. Sections \ref{Section-delta-to-0}, \ref{Section-Epsilon} and \ref{Section-limit-n}  deal with  several  limit passages w.r.t. approximation parameters.    Finally, we conclude in Section \ref{Section-Conclusion} with some remarks and we give some auxiliary results in Appendix \ref{Section-Appendix}.




\section{Concept of dissipative solution and main result}\label{Sec-weak-formulation}

Let us introduce the concept of {\em dissipative solutions}  by following the work of Feireisl et al. \cite{Dissipative-Feireisl}, which satisfy the system in the distributional
sense but with an extra ``turbulent'' term $\Re$ in the balance of momentum \eqref{momen-eq-1} that we may
call ``Reynolds stress''.

\vspace*{.15cm}
Throughout the paper, $C$ denotes a generic positive constant that does not depend on time $T$ or any underlying approximation parameters, but may depend on the given data and model parameters. If needed, we shall specify the actual  dependence of a constant on certain quantities.    By $\C_{\weak}([0,T]; X)$, we denote the  space which contains the functions which are continuous w.r.t. weak topology in the Banach space $X.$  We also use the following notations -- for any two matrices $A, B\in \R^{3\times 3}$,
\begin{align*}
A:B = \sum_{1\leq i,j\leq 3} A_{ij} B_{ij}, \ \ \ \nabla_x A \odot \nabla_x B =  (\partial_i A : \partial_j B)_{1\leq i,j \leq 3} .
\end{align*}
The inner product between $A, B$ in $L^2(\mathcal O)$ (for any open set $\mathcal O\subseteq \R^3$) is defined by 
\begin{align*}
 \int_{\mathcal{O}} A : B  = \int_{\mathcal O} \tr(A B)   .
\end{align*}
 Finally, by $\mathcal M^+(\overline{\Omega})$, we denote the set of all positive Borel measures on $\overline{\Omega}$, and if the symbol is  $\mathcal M^+(\overline{\Omega}; \R^{3\times 3})$, the measures are tensor-valued.

\subsection{Pressure--density equation of state}\label{Section-Pressure}
For the barotropic compressible flows,  the most commonly used  (isentropic) pressure--density relation is: 
\begin{equation}\label{S1}
p(\rho) = a \rho^\gamma, \ a > 0, \ \gamma > 1 \ \mbox{with associated pressure potential}\ P(\rho) = \frac{a}{\gamma - 1} \rho^\gamma.
\end{equation}
But,  in the present article, we  choose more general pressure--density equation of state as addressed in 
\cite{Dissipative-Feireisl};  more specifically,  
\begin{align} \label{S2}
\begin{dcases} 
p \in 
\C[0, \infty) \cap \C^2(0, \infty),\ 
p(0) = 0, \ p'(\rho) > 0 \ \mbox{for}\ \rho > 0, \\
\mbox{the pressure potential}\ P 
\ \mbox{determined by}\ P'(\rho) \rho - P(\rho) = p(\rho)\ \mbox{satisfies}\
P(0) = 0, \\ 
\mbox{and}\ P - \underline{a} p,\ \overline{a} p - P\ \mbox{are convex functions 
for certain constants}\ \underline{a} > 0, \ \overline{a} > 0.
\end{dcases} 
\end{align} 
In this case,
the  potential $P$ satisfying \eqref{S2} possesses some coercivity property 
similar to \eqref{S1}. Specifically, one can show that 
\begin{equation} \label{S2a}
P(\rho) \geq \tilde{a} \rho^\gamma \ \, \mbox{for certain} \ \,  \tilde{a} > 0,\ \, \gamma > 1 
\, \ \mbox{and all}\  \, \rho \geq 1.
\end{equation}
Indeed, as $\overline{a}p - P$ is a convex function and $P$ is strictly convex, one can obtain that
\begin{align*}
 p''(\rho) \geq \frac{P''(\rho)}{\overline{a}} = \frac{p'(\rho)}{\overline{a} \rho}, \ \  \  \forall  \rho > 0 , 
\end{align*}
and hence, the pressure is convex, cf.  \cite[Section 2.1.1]{Dissipative-Feireisl}. We also point out that,  due to its $\C^2(0,\infty)$-regularity, the pressure
is locally Lipschitz-continuous in $(0,\infty)$.

 Finally,  observe that we can treat the case $\gamma \in (1, \frac{3}{2}]$ whilst dealing with  dissipative solutions. But for the usual notion of weak solutions for compressible fluids, the natural restriction is  $\gamma>\frac{3}{2}$ for pressure law \eqref{S1}  
-- we refer for instance the books \cite{Feireisl-Novotny,Straskava-Antonin} for barotropic and non-barotropic compressible fluids, and for compressible active nematic liquid crystals, we refer \cite{C-Majumdar-W-Z-Comp}.

\subsection{Continuity equation -- mass conservation}
The density $\rho \geq 0$ a.e. in $(0,T)\times \Omega$, 
\begin{align*}
\rho \in \C_{\textnormal{weak}}([0,T];  L^\gamma(\Omega)) \cap L^{\gamma}(0,T; L^\gamma( \partial \Omega; |\uvect_B \cdot \mathbf n| \d S_x )), 
\end{align*}
the momentum  satisfies 
\begin{align*}
\rho \uvect \in L^\infty(0,T; L^{\frac{2\gamma}{\gamma+1}}(\Omega;\R^3)  ) ,
\end{align*}
and the following integral identity holds 
\begin{align}\label{weak-conti}
\int_\Omega \rho(\tau) \phi(\tau) \d x 
+ \int_0^\tau \int_{\Gamma_{\out}} \phi \rho \uvect_B \cdot \mathbf n \d S_x \dt + \int_0^\tau \int_{\Gamma_{\inn}} \phi \rho_B \uvect_B \cdot \mathbf n \d S_x \dt 
\notag \\ 
= \int_{\Omega} \rho_0 \phi(0,\cdot) + \int_0^\tau \int_\Omega \left[\rho \partial_t \phi + \rho \uvect \cdot \nabla_x \phi  \right] \d x \d t 
\end{align} 
for any $\tau \in [0,T]$ and $\phi \in \C^1([0,T]\times \overline{\Omega})$, with $\rho(0, \cdot)= \rho_0$ in $\Omega$, and $\rho_B \in \C^1(\mathbb R)$, $\uvect_B \in \C^1(\mathbb R^3;\R^3)$.

\subsection{Momentum balance} 

Recall that $\mathcal{M}^+(\overline{\Omega}; \R^{3\times 3}_{\sym})$ denotes the set of all positive semi-definite tensor valued measures on $\overline{\Omega}$. There exist 
\begin{align*}
\Ss \in L^1((0,T)\times \Omega; \R^{3\times 3}_{\sym}), \ \ \Re \in L^\infty(0,T; \mathcal M^+(\overline{\Omega}; \R^{3\times 3}_{\sym}) )
\end{align*}
such that the integral identity
\begin{align}\label{weak-momentum}
\int_\Omega (\rho \uvect) (\tau) \cdot \boldvphi(\tau) \d x =
\int_0^\tau \int_\Omega \left[\rho \uvect \cdot \partial_t \boldvphi + \rho \uvect \otimes \uvect : \nabla_x \boldvphi + p(\rho) \Div_x \boldvphi  - \Ss : \nabla_x \boldvphi \right] \d x \dt \notag \\
+ \int_0^\tau \int_\Omega \nabla_x \boldvphi : \d \Re (t) \d t - \int_0^\tau \int_\Omega \left[G(Q)\mathbb I_3 - \nabla_x Q \odot \nabla_x Q \right] : \nabla_x \boldvphi \d x \d t  \notag \\
- \int_0^\tau \int_\Omega \left[ Q \Delta_x Q - \Delta_x Q Q + {\sigma_* c^2 Q }  \right] : \nabla_x \boldvphi \d x \d t
+ \int_\Omega (\rho \uvect)_0 \cdot \boldvphi(0) \dx 
\end{align} 
holds for any $\tau \in [0,T]$ and any test function $\boldvphi \in \C^1_c([0,T]\times  \Omega; \R^3)$,   where $(\rho \uvect)(0,\cdot)=  (\rho \uvect)_0$ in $\Omega$.

\subsection{Concentration of active particles} 
The following integral identity holds 
\begin{align}\label{weak-concen}
\int_\Omega c(\tau)\psi (\tau) \dx - \int_0^\tau \int_\Omega c \partial_t \psi \dx\dt + \int_0^\tau \int_\Omega (\uvect \cdot \nabla_x) c \psi \dx\dt \notag 
\\ = \int_\Omega c_0 \psi(0) \dx - D_0 \int_0^\tau \int_\Omega \nabla_x c \cdot \nabla_x \psi \dx\dt 
\end{align}
for any $\tau \in [0,T]$ and any test function $\psi \in \C^1([0,T]\times \overline \Omega)$, with $c(0,\cdot)=c_0$ in $\Omega$. 

\subsection{Nematic tensor order parameter} The nematic tensor order parameter $Q$ satisfies 
\begin{align}\label{weak-nematic-tensor}
&\int_\Omega Q(\tau):\Psi (\tau) \dx - \int_0^\tau \int_\Omega Q : \partial_t \Psi \dx\dt + \int_0^\tau \int_\Omega (\uvect \cdot \nabla_x) Q : \Psi \dx\dt  \notag \\
&\quad + \int_0^\tau \int_\Omega [Q\Lambda - \Lambda Q]: \Psi \dx\dt 
= \int_\Omega Q_0 : \Psi(0) \dx - \Gamma \int_0^\tau \int_\Omega H[Q,c] :  \Psi \dx\dt 
\end{align}
for any $\tau \in [0,T]$ and any test function $\Psi \in \C^1_c([0,T]\times  \Omega; \R^{3\times 3})$, with $Q(0,\cdot)=Q_0$ in $\Omega$. Moreover, it holds that $Q\in S^3_0$.

\subsection{Energy inequality and defect compatibility condition} 

We consider the extended data $\uvect_B \in \C^1(\R^3; \R^3)$ and $\uvect-\uvect_B\in L^{\sfrac{4}{3}}(0,T; W^{1, \sfrac{4}{3}}_0(\Omega; \R^3))$. We also introduce the energy defect measure 
\begin{align}
\mathfrak{E} \in L^\infty(0,T; \mathcal{M}^+ (\overline{\Omega}))  .
\end{align}
The energy inequality  reads as 
\begin{align}\label{Energy-inequality-main}
&\frac{\d}{\d t} \int_{\Omega} \left(\frac{1}{2}\rho |\uvect - \uvect_B|^2 + P(\rho) + \frac{1}{2} |c|^2 + \frac{1}{2}|Q|^2 + \frac{1}{2}|\nabla_x Q|^2 + \frac{c_*}{4} |Q|^4 \right)   \d x  + \frac{\d}{\d t}\int_{\Omega} \mathfrak{E} \d x  \notag \\
& + \frac{1}{2}\int_{\Omega} \left[F(\D_x\uvect) + F^*(\Ss) \right]\d x  
+  \frac{D_0}{2} \int_\Omega |\nabla_x c|^2 \d x \notag + \frac{\Gamma}{4} \int_{\Omega} |\Delta_x Q|^2  \d x  + 
\frac{c_*^2 \Gamma}{2} \int_{\Omega} |Q|^6\dx    \notag \\
& + \int_{\Gamma_{\out}} P(\rho) \uvect_B \cdot \mathbf n  \d S_x  + \int_{\Gamma_{\inn}} P(\rho_B) \uvect_B \cdot \mathbf n \d S_x  
\notag \\
  \leq  & C \int_{\Omega} \left(\frac{1}{2} |c|^2 + \frac{1}{2}|Q|^2 + \frac{1}{2}|\nabla_x Q|^2 + \frac{c_*}{4} |Q|^4 \right)   \d x  + \int_\Omega \Ss : \nabla_x \uvect_B \d x  \notag \\ 
  & - \int_{\Omega} \left[ \rho \uvect \otimes \uvect + p(\rho) \mathbb{I}_3\right] : \nabla_x \uvect_B \d x 
  + \int_\Omega \rho \uvect \cdot (\uvect_B \cdot \nabla_x \uvect_B ) \d x - 
  \int_{\Omega} \nabla_x \uvect_B : \d \Re  \notag \\
 & - \frac{1}{2} \int_{\partial \Omega} \left(\frac{1}{2} |Q_B|^2 +  \frac{c_*}{4} \tr^2(Q_B^2) \right) \uvect_B \cdot \mathbf n   \d S_x  
  + 2c_* \Gamma \int_{\partial \Omega} |Q_B|^2 Q_B : \nabla_x Q_B \cdot \mathbf n \d S_x + C ,
\end{align}
for some constant $C>0$ that depends at most on the given initial and boundary data (see Section \ref{Section-Apriori} for the derivation of  energy inequality  at the a-priori level for regular enough solutions).

Finally, we impose a compatibility condition between the energy defect $\mathfrak E$ and Reynolds stress $\Re$, specifically
\begin{align}
\underline{d} \mathfrak{E} \leq \tr[\Re] \leq \overline{d} \mathfrak{E} 
\end{align}
for certain constants $0<\underline{d}\leq \overline{d}$.  See Section \ref{Section-subsub-Momentum}  about why this condition is appearing in our analysis.



\begin{defi}[\textbf{Dissipative solution}]\label{Definition-sol}
Let $\Omega \subset \mathbb R^3$ be a bounded Lipschitz domain. Then, $(\rho, \uvect, c, Q)$ is said to be {dissipative solution} of the system \eqref{dens-eq-1}--\eqref{nematic-tensor-eq-1} with prescribed initial and boundary conditions if:  
\begin{itemize}
\item $\rho \geq 0$   a.e. in  $(0,T)\times \Omega$, and 
\begin{align*}
\ \rho \in \C_{\textnormal{weak}}([0,T];  L^\gamma(\Omega)) \cap L^{\gamma}(0,T; L^\gamma( \partial \Omega; |\uvect_B \cdot \mathbf n| \d S_x )), \ \gamma>1. 
\end{align*}
\begin{align*}
\uvect \in L^{\sfrac{4}{3}}(0,T; W^{1, \sfrac{4}{3} }(\Omega; \R^3)) , \  (\uvect- \uvect_B) \in L^{\sfrac{4}{3}}(0,T; W^{1, \sfrac{4}{3} }_0(\Omega; \R^3)), 
\end{align*}
\begin{align*}
 \rho \uvect \in \C_{\textnormal{weak}}([0,T]; L^{\frac{2\gamma}{\gamma+1}}(\Omega ; \R^3)).
\end{align*}

\item $c>0$, $c\in L^\infty((0,T)\times \Omega) \cap L^2(0,T; H^1(\Omega ))$, and $Q\in L^\infty(0,T; H^1(\Omega) ) \cap L^2(0,T; H^2(\Omega ))$, $Q\in S^3_0$ a.e. in $(0,T)\times \Omega$.

\item There exist 
\begin{align*}
\Ss \in L^1((0,T)\times \Omega; \R^{3\times 3}_{\sym}), \ \Re \in L^\infty(0,T; \mathcal M^+(\overline{\Omega}; \R^{3\times 3}_{\sym}) ) , \ \mathfrak{E} \in L^\infty(0,T; \mathcal M^+(\overline{\Omega})),
\end{align*}
such that the relations \eqref{weak-conti}--\eqref{weak-nematic-tensor} are satisfied for any $0\leq \tau \leq T$,  and the energy inequality \eqref{Energy-inequality-main} holds. 

\end{itemize}

\end{defi}

Next, we state the main result of the present article. 

\begin{theorem}[\textbf{Global existence of dissipative solution}]\label{Theorem-Main}
Let $T>0$ and $\Omega\subset \mathbb R^3$ be a bounded Lipschitz domain.  Suppose that $p=p(\rho)$ is given by \eqref{S2}  and  $F=F(\mathbb D_x \uvect)$ is given by \eqref{F-convex}--\eqref{Coercivity-F}. 
Then the compressible active nematic liquid crystal system \eqref{dens-eq-1}--\eqref{nematic-tensor-eq-1},  with initial conditions \eqref{initial-concen}--\eqref{initial-Q} and boundary conditions \eqref{inflow-outflow}--\eqref{boundary-c}, admits at least one dissipative  solution in the sense of Definition \ref{Definition-sol}. 
\end{theorem}

\subsection*{Strategy of the proof}
The proof is based on a three-level approximation scheme.
\begin{itemize}
\item[(i)] Artificial viscosity $\veps \nabla_x \rho$ for $\veps >0$: we add this term to the continuity equation \eqref{dens-eq-1} which gives parabolic regularization effect. Accordingly, in the momentum equation \eqref{momen-eq-1} we need to add the term $\veps \nabla_x \rho \cdot \nabla_x \uvect$. 

\item[(ii)] Regularization the convex potential $F$ by  $F_\delta$  for $\delta>0$: see \eqref{F-delta}.

\item[(iii)] Finite-dimensional approximation ($n$-layer): the momentum equation \eqref{momen-eq-1} has to be solved by means of a Galerkin approximation.

\end{itemize}

Then, after obtaining the approximate solutions for  fixed $\delta, \veps, n$, we pass to the limit w.r.t. the aforementioned parameters. First, we analyse the limiting behaviour  w.r.t.  $\delta\to 0$ which is followed by the limit passage of $\veps \to 0$. In the final level, we pass to the limit w.r.t. $n\to \infty$ to prove the existence of dissipative solutions for our compressible active NLC model.

\section{A-priori estimates}\label{Section-Apriori}
This section is devoted to a-priori estimates which provide us with a proper  energy inequality for our compressible system. By formally multiplying the equation \eqref{momen-eq-1} by $\uvect-\uvect_B$, \eqref{concen-eq-1} by $c$ and \eqref{nematic-tensor-eq-1} by $-(\Delta_x Q - Q - c_* Q \tr(Q^2))$, and adding, we get 
\begin{align}\label{Energy-ineq}
&\frac{\d}{\d t} \int_{\Omega} \left(\frac{1}{2}\rho |\uvect - \uvect_B|^2 + P(\rho) + \frac{1}{2} |c|^2 + \frac{1}{2}|Q|^2 + \frac{1}{2}|\nabla_x Q|^2 + \frac{c_*}{4} |Q|^4 \right)   \d x  \notag \\
& + \int_{\Omega} \left[F(\D_x\uvect) + F^*(\Ss) \right]\d x  
+  D_0 \int_\Omega |\nabla_x c|^2 \d x \notag + \Gamma \int_{\Omega} |\Delta_x Q|^2  \d x  + 
c_*^2 \Gamma \int_{\Omega} |Q|^6\dx    \notag \\
& + \int_{\Gamma_{\out}} P(\rho) \uvect_B \cdot \mathbf n  \d S_x  + \int_{\Gamma_{\inn}} P(\rho_B) \uvect_B \cdot \mathbf n \d S_x  \notag  \\
 = &  - \int_{\Omega} [\rho \uvect \otimes \uvect + p(\rho) \mathds{I}_3] : \nabla_x \uvect_B \d x   + \int_{\Omega} \rho \uvect \cdot (\uvect_B \cdot \nabla_x \uvect_B) \d x  + \int_\Omega \Ss : \nabla_x \uvect_B \dx  \notag \\
 &  - \int_{\Omega} G(Q)\mathds{I}_3 : \nabla_x \uvect  \d x  + \int_{\Omega} G(Q)\mathds{I}_3 : \nabla_x \uvect_B  \d x     -  \int_{\Omega} \Div_x (\nabla_x Q \odot \nabla_x Q)   \cdot \uvect \d x   \notag \\
 &  + \int_{\Omega} \Div_x (\nabla_x Q \odot \nabla_x Q)   \cdot \uvect_B \d x
 + \int_{\Omega} \Div_x (Q\Delta_x Q - \Delta_x Q Q ) \cdot \uvect \d x  \notag \\
 &-   \int_{\Omega} \Div_x (Q\Delta_x Q - \Delta_x Q Q ) \cdot \uvect_B \d x 
- \int_{\Omega} (\uvect \cdot \nabla_x c ) \, c  \d x  \notag \\
& + \int_{\Omega} (Q \Lambda - \Lambda Q) : \Delta_x Q \d x 
  - \int_{\Omega} (Q \Lambda - \Lambda Q) : ( Q + c_* Q |Q|^2)\d x   \notag \\
& + \int_{\Omega} (\uvect \cdot \nabla_x Q  ) :  \Delta_x Q \d x 
- \int_{\Omega} (\uvect \cdot \nabla_x Q ) : (Q + c_* Q |Q|^2) \d x \notag \\
  & + \Gamma \int_{\Omega} \frac{c-c_*}{2} Q : \left(\Delta_x Q - Q - c_* Q |Q|^2 \right) \d x - b \int_{\Omega} \left(Q^2-\frac{|Q|^2}{3}\mathbb I_3\right) : \left(\Delta_x Q - Q - c_* Q|Q|^2 \right) \d x \notag \\
&   + \Gamma \int_{\Omega} \Delta_x Q Q \d x   + 2c_*\Gamma \int_\Omega Q|Q|^2  : \Delta_x Q \d x  - c_* \int_\Omega |Q|^4 \d x  \notag \\
&{-\sigma_* \int_{\Omega} c^2 Q : \nabla_x \uvect \d x +  \sigma_* \int_{\Omega} c^2 Q : \nabla_x \uvect_B} \d x  \notag \\
& = \sum_{j=1}^{21} I_{j} .
\end{align}

We first calculate $$I_8+  I_{9}+ I_{11}.$$ Using the fact that $\uvect|_{\partial \Omega}=\uvect_B$ and Lemma \ref{lema-aux-1}, we deduce that
\begin{align}\label{E-9-10-12} 
&I_8 +  I_{9}+ I_{11} \notag \\
= & \int_{\Omega} \Div_x (Q\Delta_x Q - \Delta_x Q Q) \cdot \uvect \d x  - \int_{\Omega} \Div_x (Q\Delta_x Q - \Delta_x Q Q) \cdot \uvect_B \d x + \int_{\Omega} (Q\Lambda - \Lambda Q) : \Delta_x Q \d x\notag \\
= &  {\int_{\Omega} (Q\Delta_x Q- \Delta_x Q Q) : \nabla_x^\top \uvect_B \d x} .  
\end{align}

Next, we compute the following,
\begin{align}
& I_{6} + I_{7} + I_{13} + I_{14}  \notag \\
= & -\int_{\Omega} \Div_x (\nabla_x Q \odot \nabla_x Q ) \cdot \uvect \d x + 
\int_{\Omega} \Div_x (\nabla_x Q \odot \nabla_x Q ) \cdot \uvect_B \d x \notag \\
& + \int_\Omega (\uvect \cdot \nabla_x Q) : \Delta_x Q \d x - \int_\Omega (\uvect\cdot \nabla_x Q) : (Q + c_* Q|Q|^2) \d x  \notag  \\
= & - \int_{\Omega} \partial_j \partial_i Q^{kl} \partial_j Q^{kl} u_i \d x - \int_{\Omega}  \partial_i Q^{kl} \partial^2_j  Q^{kl} u_i   \d x
+ \int_{\Omega} \Div_x (\nabla_x Q \odot \nabla_x Q ) \cdot \uvect_B \d x \notag \\
&
+ \int_{\Omega} u_i \partial_i Q^{kl} \Delta_x Q^{kl} \d x
 - \int_{\Omega}  u_i \partial_i Q^{kl} (Q^{kl} + c_* Q^{kl} |Q|^2 ) \d x \notag \\
=& - \frac{1}{2}\int_{\Omega} \partial_i |\partial_j Q^{kl}|^2 u_i \d x - \int_\Omega u_i \partial_i \left(\frac{1}{2}|Q|^2 + \frac{c_*}{4} \tr^2(Q^2)  \right) \d x \notag \\
 & + \frac{1}{2} \int_\Omega \partial_i |\partial_j Q^{kl}|^2 u_{B,i} \d x + \int_{\Omega} \partial_i Q^{kl} \partial^2_j  Q^{kl} u_{B,i} \d x \notag \\
 =& \int_{\Omega} \left(\frac{1}{2} |\nabla_x Q|^2 + \frac{1}{2} |Q|^2 + \frac{c_*}{4} \tr^2(Q^2) \right) \Div_x \uvect   \d x   - \frac{1}{2} \int_{\Omega} |\nabla_x Q|^2 \Div_x \uvect_B \d x \notag \\
 & + \int_\Omega (\uvect_B \cdot \nabla_x Q)   :  \Delta_x Q \d x  - {\int_{\partial \Omega} \left(\frac{1}{2}|Q_B|^2 + \frac{c_*}{4}\tr^2(Q_B^2)  \right) \uvect_B \cdot \mathbf n \d S_x }.
\end{align}
But 
$$I_4 = - \int_{\Omega} G(Q) \mathbb{I}_3 : \nabla_x \uvect \d x = -\int_{\Omega} \left(\frac{1}{2} |\nabla_x Q|^2 + \frac{1}{2} |Q|^2 + \frac{c_*}{4} \tr^2(Q^2) \right) \Div_x \uvect   \d x , $$
and thus, we have 
\begin{align}\label{E-5-15}
& I_4 + I_6 + I_{7} + I_{13} + I_{14} \notag \\
= & - \frac{1}{2} \int_\Omega |\nabla_x Q|^2 \Div_x \uvect_B \d x + \int_{\Omega} (\uvect_B \cdot \nabla_x Q ) : \Delta_x Q \d x  \notag \\
& - \int_{\partial \Omega} \left(\frac{1}{2} |Q_B|^2 + \frac{c_*}{4} \tr^2(Q_B^2) \right)  \uvect_B \cdot \mathbf n \d S_x .   
\end{align}

Next, observe that 
\begin{align}\label{E-13}
I_{12} & =  -\int_\Omega ( Q \Lambda - \Lambda Q) : (Q + c_* Q |Q|^2) \d x \notag \\
& = - \int_\Omega (Q \Lambda + \Lambda Q): (Q + c_* Q |Q|^2) \d x + 2 \int_\Omega \Lambda Q : (Q + c_* Q |Q|^2)  \d x = 0,
\end{align}
where we have used the fact that $Q$ is symmetric and $\Lambda$ is skew-symmetric matrices.

We further find that 
\begin{align*} 
I_{10}= -\int_{\Omega} (\uvect \cdot \nabla_x c) c \d x  = -\frac{1}{2}\int_{\Omega} \uvect \cdot \nabla_x |c|^2 = \frac{1}{2}\int_{\Omega} |c|^2 \Div_x \uvect \d x  - \frac{1}{2}\int_{\partial \Omega} |c|^2 \uvect_B \cdot \mathbf n \d S_x .
\end{align*}
Therefore,
\begin{align}\label{E-11}
|I_{10}| &\leq \frac{1}{2} \int_{\Omega}|c|^2 \Div_x \uvect \d x + C\|\uvect_B\|_{L^\infty(\partial\Omega)} \|c\|^2_{L^2(\partial \Omega)} \notag   \\
& \leq  \frac{1}{2} \int_{\Omega}|c|^2 \Div_x \uvect \d x + \frac{D_0}{2}\|\nabla_x c\|^2_{L^2(\Omega)} + C(D_0, \uvect_B) \|c\|^2_{L^2(\Omega)} ,
\end{align}
where we have used the trace inequality reported in Lemma \ref{Trace-them}. 

Now, we find the following estimates. First, note that $c \in L^\infty((0,T)\times \Omega)$ as follows from the maximum principle for the equation of concentration $c$, and the choice of $c_0 \in L^\infty(\Omega)$. We have 
\begin{align}
 |I_{15}| &= \left|\Gamma \int_{\Omega} \frac{(c-c_*)}{2} Q : \left(\Delta_x Q - Q - c_* Q |Q|^2\right) \d x  \right| \leq \frac{\Gamma}{6} \|\Delta_x Q\|^2_{L^2} + C\|Q\|^2_{L^2} + C \|Q\|^4_{L^4} ,  \label{E-16} \\
 |I_{16}| &= \left| b \int_\Omega \left(Q^2 - \frac{|Q|^2}{3}\mathbb I_3\right) : \left(\Delta_x Q - Q - c_* Q |Q|^2\right) \d x \right|    \notag  \\
& \leq \frac{\Gamma}{6} \|\Delta_x Q \|^2_{L^2} + C \|Q\|^4_{L^4} + C\left|  \int_{\Omega} \tr(Q^3) \d x \right| + C c_* \left| \int_\Omega \tr(Q^3) |Q|^2 \d x \right| \notag  \\
& \leq \frac{\Gamma}{6} \|\Delta_x Q \|^2_{L^2} + C \|Q\|^4_{L^4} + \frac{c^2_*\Gamma}{2}\|Q\|^6_{L^6} + C\|Q\|^2_{L^2} ,  \label{E-17} \\
 |I_{17}| &= \left|\Gamma \int_{\Omega} \Delta_x Q Q \d x  \right| \leq \frac{\Gamma}{6} \int_{\Omega} |\Delta_x Q|^2 \d x + C \|Q\|^2_{L^2}   \label{E-18} . 
 \end{align}
 
We also obtain that
\begin{align*}
&I_{18}  = 2 c_* \Gamma \int_\Omega Q|Q|^2 : \Delta_x Q \d x = 2 c_* \Gamma \int_\Omega \partial^2_{j} Q^{kl}  Q^{kl} |Q|^2 \d x   \notag \\
& = - 2 c _* \Gamma \int_\Omega |\partial_j Q^{kl}|^2 |Q|^2 \d x - 2 c_* \Gamma \int_\Omega \partial_j Q^{kl} Q^{kl} \partial_j (\tr(Q^2)) +  2c_* \Gamma \int_{\partial \Omega} |Q_B|^2 Q_B : \nabla_x Q_B \cdot \mathbf n  \d S_x \\
& = - 2 c _* \Gamma \int_\Omega |\nabla_x Q|^2 |Q|^2 \d x - 2 c_* \Gamma \int_\Omega \partial_j Q^{kl}  Q^{kl} \partial_j (\tr(Q^2)) +  2c_* \Gamma \int_{\partial \Omega} |Q_B|^2 Q_B : \nabla_x Q_B \cdot \mathbf n \d S_x  \\
& = \underbrace{ - 2c_* \Gamma \int_\Omega |\nabla_x Q|^2 |Q|^2 \d x - c_* \Gamma \int_\Omega |\nabla_x (\tr(Q^2))|^2 \d x }_{I^1_{18}} +  \underbrace{ 2 c_* \Gamma \int_{\partial \Omega} |Q_B|^2 Q_B : \nabla_x Q_B \cdot \mathbf n  \d S_x}_{I^{2}_{18}}.  
\end{align*}
Observe that 
\begin{align}\label{E-19-1}
I^{1}_{18} \leq 0,
\end{align}
and the boundary term $I^2_{18}$ is bounded:
\begin{align}\label{E-19-2}
|I^2_{18}| \leq C \left(\|Q_B\|^6_{L^6(\partial \Omega)} + \|\nabla_x Q_B\|^2_{L^2(\partial \Omega}\right) ,
\end{align}
since $Q_B \in H^{\sfrac{3}{2}}(\partial \Omega)$.

Combining the estimates \eqref{E-9-10-12}--\eqref{E-19-1}, the energy inequality \eqref{Energy-ineq} reduces to
\begin{align}\label{Energy-ineq-2}
&\frac{\d}{\d t} \int_{\Omega} \left(\frac{1}{2}\rho |\uvect - \uvect_B|^2 + P(\rho) + \frac{1}{2} |c|^2 + \frac{1}{2}|Q|^2 + \frac{1}{2}|\nabla_x Q|^2 + \frac{c_*}{4} |Q|^4 \right)   \d x  \notag \\
& + \int_{\Omega} \left[F(\D_x\uvect) + F^*(\Ss) \right]\d x  
+  \frac{D_0}{2} \int_\Omega |\nabla_x c|^2 \d x \notag + \frac{\Gamma}{2} \int_{\Omega} |\Delta_x Q|^2  \d x  + 
\frac{c_*^2 \Gamma}{2} \int_{\Omega} |Q|^6\dx    \notag \\
& + \int_{\Gamma_{\out}} P(\rho) \uvect_B \cdot \mathbf n  \d S_x  + \int_{\Gamma_{\inn}} P(\rho_B) \uvect_B \cdot \mathbf n \d S_x  \notag  \\
 \leq &  - \int_{\Omega} [\rho \uvect \otimes \uvect + p(\rho) \mathds{I}_3] : \nabla_x \uvect_B \d x   + \int_{\Omega} \rho \uvect \cdot (\uvect_B \cdot \nabla_x \uvect_B) \d x  + \int_\Omega \Ss : \nabla_x \uvect_B \dx  \notag \\
 & + \int_{\Omega} G(Q) \mathbb{I}_3 : \nabla_x \uvect_B \d x - \frac{1}{2}\int_\Omega  |\nabla_x Q|^2 \Div_x \uvect_B \d x + \int_\Omega (\uvect_B \cdot \nabla_x Q) : \Delta_x Q \d x \notag \\
 & +\frac{1}{2}\int_{\Omega}|c|^2 \Div_x \uvect \d x 
  + \int_\Omega (Q \Delta_x Q - \Delta_x Q Q): \nabla^\top_x \uvect_B \d  x   \notag \\
  & -\sigma_* \int_{\Omega} c^2 Q : \nabla_x \uvect \d x +  \sigma_* \int_{\Omega} c^2 Q : \nabla_x \uvect_B \d x \notag \\  
&  - \frac{1}{2} \int_{\partial \Omega} \left(\frac{1}{2} |Q_B|^2 +  \frac{c_*}{4} \tr^2(Q_B^2) \right) \uvect_B \cdot \mathbf n   \d S_x + 2c_* \Gamma \int_{\partial \Omega} |Q_B|^2 Q_B : \nabla_x Q_B \cdot \mathbf n \d S_x \notag \\ 
  &   + C\left(\|c\|^2_{L^2(\Omega)}+  \|Q\|^2_{L^2(\Omega)} + \|Q\|^4_{L^4(\Omega)} \right) ,
\end{align}
where the constant $C>0$ depends at most on the system parameters and given  data. 

%

\vspace*{.2cm}

A few further computations are needed to conclude the derivation of the energy inequality in \eqref{Energy-ineq-2}.

\begin{itemize}

\item[(i)] We begin with the following:
\begin{align}\label{esti-1}
\int_\Omega |G(Q)| |\nabla_x \uvect_B|\leq C\|\nabla_x \uvect_B \|_{L^\infty(\Omega)}\int_\Omega \left(|\nabla_x Q|^2 + |Q|^2 + |Q|^4 \right)\d x  ,
\end{align}
and 
\begin{align}
\frac{1}{2}\left| \int_{\Omega} |\nabla_x Q|^2 \Div_x \uvect_B \d x \right| \leq C \|\Div_x \uvect_B \|_{L^\infty(\Omega)} \int_{\Omega} |\nabla_x Q|^2 \d x.
\end{align}

\item[(ii)] Further, we calculate 
\begin{align}\label{esti-3}
\int_\Omega| (\uvect_B \cdot \nabla_x Q) : \Delta_x Q | \leq   \frac{\Gamma}{8} \int_\Omega |\Delta_x Q|^2 + C(\Gamma)\|\uvect_B\|_{L^\infty(\Omega)} \int_\Omega |\nabla_x Q|^2  ,
\end{align}
and 
\begin{align}\label{esti-4}
\left|\int_{\Omega} (Q\Delta_x Q - \Delta_x Q Q): \nabla^\top_x \uvect_B \right| &\leq 2\|\nabla_x \uvect_B\|_{L^\infty(\Omega)} \int_{\Omega} |Q | |\Delta_x Q| \notag \\
&\leq \frac{\Gamma}{8} \|\Delta_x Q\|^2_{L^2(\Omega)} + C(\uvect_B, \Gamma) \|Q\|^2_{L^2(\Omega)} . 
\end{align}

\item[(iii)] Using  the 
$L^q$-version of trace-free Korn's inequality given in Lemma \ref{Lemma-Korn-trace-free} and the
convexity \eqref{F-convex}--\eqref{Coercivity-F},  one has 
\begin{align}\label{Korn-trace-free-0}
\| \nabla_x \uvect \|^{\sfrac{4}{3}}_{L^{\sfrac{4}{3}}(\Omega)} \leq \kappa_0 \left\|\D_x \uvect - \frac{1}{3} \Div_x \uvect \mathbb{I}_3 \right\|^{\sfrac{4}{3}}_{L^{\sfrac{4}{3}}(\Omega)} +  C(\uvect_B) \\
\label{Korn-trace-free-1}
\leq \kappa_1 \int_{\Omega} \left[ F (\mathbb D_x \uvect) + F^*(\Ss) \right] \d x + C(\uvect_B), 
\end{align}
for some positive constants $\kappa_0$, $\kappa_1$.   Observe that, in the inequality \eqref{Korn-trace-free-0}, there is an extra constant  term $C(\uvect_B)$ depending on $\uvect_B$ (unlike the Korn's inequality \eqref{Trace-free-Korn-ineq}), and this is due to the fact that we consider nonhomogeneous Dirichlet data for the velocity $\uvect$.

With this,  we compute the following: 
\begin{align}\label{esti-6}
\left|\frac{1}{2}\int_{\Omega} |c|^2 \Div_x \uvect  \d x \right| &\leq C \|\nabla_x \uvect\|_{L^{\sfrac{4}{3}}(\Omega)} \|c\|^2_{L^8(\Omega)} \notag \\
& \leq   \frac{1}{4\kappa_1} \|\nabla_x \uvect\|^{\sfrac{4}{3}}_{L^{\sfrac{4}{3}}(\Omega)} + C(\kappa_1) \|c\|^8_{L^8(\Omega)} \notag \\
& \leq \frac{1}{4}\int_{\Omega} \left[F (\mathbb D_x \uvect) + F^*(\Ss) \right] \d x + C(\uvect_B, \kappa_1) (1+\|c\|^8_{L^\infty(\Omega)}), 
\end{align}
where we have used the inequality \eqref{Korn-trace-free-1} and the fact that $c\in L^\infty((0,T)\times \Omega)$.

\item[(iv)] In terms of Cauchy-Schwarz inequality and \eqref{Korn-trace-free-1}, we now estimate the following:
\begin{align}\label{esti-7}
\sigma_*\left|\int_{\Omega} c^2 Q: \nabla_x \uvect \right| &\leq \sigma_*\|c\|^2_{L^\infty(\Omega)}\int_{\Omega} |Q||\nabla_x \uvect| \leq \frac{1}{4\kappa_1} \|\nabla_x \uvect\|^{\sfrac{4}{3}}_{L^{\sfrac{4}{3}}(\Omega)} + C(\sigma_*, \kappa_1)  \|Q\|^4_{L^4(\Omega)} \notag 
\\
& \quad \leq \frac{1}{4}\int_{\Omega} \left[ F (\mathbb D_x \uvect) + F^*(\Ss) \right] \d x + C(\uvect_B, \sigma_*, \kappa_1)(1+\|Q\|^4_{L^4(\Omega)} ) . 
\end{align}
Moreover, we see that
\begin{align}\label{esti-8}
\sigma_*\left|\int_{\Omega} c^2 Q: \nabla_x \uvect_B \right| \leq \sigma_* \|c\|^2_{L^\infty(\Omega)}\left( \|Q\|^2_{L^2(\Omega)} + \|\nabla_x \uvect_B\|^2_{L^2(\Omega)}\right)   .
\end{align}

\end{itemize}

 Finally, after combining the estimates \eqref{esti-1}--\eqref{esti-8} in \eqref{Energy-ineq-2}, we have 
\begin{align}\label{Energy-ineq-4}
&\frac{\d}{\d t} \int_{\Omega} \left(\frac{1}{2}\rho |\uvect - \uvect_B|^2 + P(\rho) + \frac{1}{2} |c|^2 + \frac{1}{2}|Q|^2 + \frac{1}{2}|\nabla_x Q|^2 + \frac{c_*}{4} |Q|^4 \right)   \d x  \notag \\
& + \frac{1}{2}\int_{\Omega} \left[F(\D_x\uvect) + F^*(\Ss) \right]\d x  
+  \frac{D_0}{2} \int_\Omega |\nabla_x c|^2 \d x \notag + \frac{\Gamma}{4} \int_{\Omega} |\Delta_x Q|^2  \d x  + 
\frac{c_*^2 \Gamma}{2} \int_{\Omega} |Q|^6\dx    \notag \\
& + \int_{\Gamma_{\out}} P(\rho) \uvect_B \cdot \mathbf n  \d S_x  + \int_{\Gamma_{\inn}} P(\rho_B) \uvect_B \cdot \mathbf n \d S_x  
\notag \\
  \leq  & C \int_{\Omega} \left( \frac{1}{2} |c|^2 + \frac{1}{2}|Q|^2 + \frac{1}{2}|\nabla_x Q|^2 + \frac{c_*}{4} |Q|^4 \right)   \d x    + \int_\Omega \Ss : \nabla_x \uvect_B \dx \notag \\ 
  & - \int_{\Omega}\left[ \rho \uvect \otimes \uvect + p(\rho) \mathbb{I}_3\right] : \nabla_x \uvect_B \d x + \int_{\Omega} \rho \uvect \cdot (\uvect_B \cdot \nabla_x \uvect_B) \d x
  \notag \\
 & - \frac{1}{2} \int_{\partial \Omega} \left(\frac{1}{2} |Q_B|^2 +  \frac{c_*}{4} \tr^2(Q_B^2) \right) \uvect_B \cdot \mathbf n   \d S_x  
  + 2c_* \Gamma \int_{\partial \Omega} |Q_B|^2 Q_B : \nabla_x Q_B \cdot \mathbf n \d S_x + C ,
\end{align}
where $C>0$ is a constant that depends on the given data   and  system parameters.  

\vspace{.5cm}

{\em The remaining sections are devoted to the proof of Theorem \ref{Theorem-Main}, which is our main result.}

\section{Existence of approximate solutions}\label{Section-Existence}

To begin with, we consider $\uvect_B \in \C^1_c(\R^3; \R^3)$, $\rho_B \in \C^1(\partial \Omega)$ and $\rho_0 \in \C^1(\R^3)$ such that $\rho_0 \equiv 0$ in $\mathbb R^3 \setminus \Omega$. We also choose $c_0\in \C^1_c(\Omega)$, $Q_0\in \C^1_c(\Omega)$ and $Q_B \in \C^1(\partial \Omega)$.


Let us consider a sequence of finite-dimensional spaces 
$X_n \subset L^2(\Omega; \mathbb R^3)$
\begin{align}\label{Space-X_n}
X_n = \span \left\{ \vc{w}_i\ \Big|\ \vc{w}_i \in \C^\infty_c(\Omega; R^d),\ i = 1,\dots, n \right\} . 
\end{align}   
Without loss of generality, we may assume that $\vc{w}_i$ are orthonormal w.r.t.  the standard scalar product in $L^2$.


\vspace*{.1cm}
 Secondly, we need to  regularize the convex potential $F$  (given by \eqref{Coercivity-F}) to make it continuously differentiable. This may be done via convolution with a family 
of regularizing kernels $\{\xi_\delta\}_{\delta>0}$ in $\R^{3\times 3}_{\sym}$, more precisely
\begin{align}\label{F-delta}
F_\delta(\mathbb{D}) = \int_{R^{3 \times 3}_{\sym}} \xi_\delta ( \mathbb{D} - \mathbb{Z}  ) F(\mathbb{Z}) \d \mathbb{Z}
- \int_{R^{3 \times 3}_{\sym}} \xi_\delta ( \mathbb{Z} ) F(\mathbb{Z}) \d \mathbb{Z} ,  \ \ \forall \, \D \in \R^{3\times 3}_{\sym} . 
\end{align}  
It is easy to check  that $F_\delta$ are convex, 
non-negative, infinitely differentiable, $F_\delta (0) = 0$. Moreover, it can be proved that $F_\delta$   satisfies 
\begin{equation}\label{F-delta-convex}
F_\delta (\mathbb{D}) \geq \mu_1 \left| \mathbb{D} - \frac{1}{d} {\rm tr}[\mathbb{D}] \mathbb{I} \right|^{\sfrac{4}{3} } - \mu_2 , 
\ \ \forall  \,  \mathbb{D} \in \mathbb R^{3\times 3}_{\sym},
\end{equation}
for some positive constants $\mu_1$ and $\mu_2$  which are  independent of $\delta \to 0$; for a detailed proof, we refer, for instance \cite[Proposition 4.2]{Basaric}.

\subsection{Artificial viscosity approximation of the equation of continuity}

In the continuity equation, we use a parabolic regularization effect, namely we consider 
\begin{align}\label{cont-eps}
\partial_t \rho + \Div_x (\rho \uvect ) = \veps \Delta_x \rho \ \mbox{ in}\ (0,T) \times \Omega,\ \veps > 0,  
\end{align}
supplemented with the boundary conditions 
\begin{equation} \label{bound-eps}
\veps \nabla_x \rho \cdot \mathbf n + (\rho_B - \rho ) [\uvect_B \cdot  \mathbf n]^- = 0 \ \mbox{ in}\ (0,T) \times \partial \Omega ,
\end{equation}
and initial condition
\begin{align}\label{initial-eps}
\rho(0, \cdot) = \rho_0 \ \text{ in } \Omega  .
\end{align}
Here, $\uvect = \vc{v} + \uvect_B$, with $\vc{v} \in \C([0,T]; X_n)$, in particular, $\uvect|_{\partial \Omega} = \uvect_B$. 
Note that for given $\vc{v}$, $\rho_B$, $\uvect_B$, this is a linear problem for the unknown $\rho$. 

As $\Omega$ is merely Lipschitz, 
we  use the weak formulation: 
\begin{equation} \label{E4}
\begin{split}
\left[ \int_\Omega  \rho \varphi  \right]_{t = 0}^{t = \tau}&= 
\int_0^\tau \int_\Omega{ \left[ \rho \partial_t \varphi + \rho \uvect \cdot \nabla_x \varphi - \veps \nabla_x \rho \cdot \nabla_x  \varphi \right] }
\dt \\ 
&- \int_0^\tau \int_{\partial \Omega} \varphi \rho \uvect_B \cdot \vc{n} \ {\rm d} S_x \dt + 
\int_0^\tau \int_{\partial \Omega} \varphi (\rho - \rho_B) [\uvect_B \cdot \vc{n}]^{-}  \ {\rm d}S_x \dt ,
\end{split}
\end{equation}
for any $0\leq \tau \leq T$ and for any test function
\[
\varphi \in L^2(0,T; W^{1,2}(\Omega)),\ \partial_t \varphi \in L^1(0,T; L^2(\Omega)).
\]
with $\rho(0,\cdot)=\rho_0$ in $\Omega$.

\vspace{.2cm}

Let us recall some important results concerning the problem \eqref{cont-eps}--\eqref{initial-eps}.

\begin{lemma}[{Crippa et al. \cite[Lemma 3.2]{Crippa-et-al}}]\label{EL1}
 Let $\Omega \subset R^3$ be a bounded Lipschitz domain.
Given $\uvect = \vc{v} + \uvect_B$, $\vc{v} \in \C([0,T]; X_n)$, the initial-boundary value problem \eqref{cont-eps}--\eqref{initial-eps} 
admits a weak solution $\rho$  in the sense of \eqref{E4}, unique in the class 
\[
\rho \in L^2(0,T; W^{1,2}(\Omega)) \cap \C([0,T]; L^2(\Omega)).
\]
The norm in the aforementioned spaces is bounded only in terms of the data $\rho_0$, $\rho_B$, $\uvect_B$ and 
\[
\sup_{t \in (0,T)}
\| \vc{v}(t, \cdot) \|_{X_n}.
\]
\end{lemma}


\begin{lemma}[{Maximum principle; Crippa et al. \cite[Lemma 3.4]{Crippa-et-al}}]\label{EL2}
Under the hypotheses of Lemma \ref{EL1}, the solution $\rho$ satisfies 
\begin{align*}
&\| \rho \|_{L^\infty((0,\tau) \times \Omega)}  \\
\leq & \max \left\{ \| \rho_0 \|_{L^\infty(\Omega)}, \| \rho_B \|_{L^\infty((0,T) \times \Gamma_{in})}, 
\| \uvect_B \|_{L^\infty((0,T) \times \partial \Omega)} \right\} \exp \left( \tau \| \Div_x \uvect \|_{L^\infty((0,\tau) \times \Omega)} \right).
\end{align*}
\end{lemma}

\vspace*{.1cm}

Next, we report another result from Abbatiello et al.  \cite[Lemma 3.3]{Dissipative-Feireisl}. 
\begin{lemma}[{Renormalization}] \label{EL3}
Under the conditions of Lemma \ref{EL2},
let $B \in \C^2(\R)$, $\chi \in W^{1, \infty}(0,T)$, and $r =\rho - \chi$.

Then, the (integrated) renormalized equation
\[
\begin{split}
\left[ \int_{\Omega}{ B(r) } \right]_{t = 0}^{t = \tau} 
= &- \int_{0}^{\tau} \int_{\Omega}{ \Div (r \uvect) B'(r) } \d t - \veps \int_{0}^{\tau} \int_{\Omega}{ |\nabla_x r|^2  B''(r) } \d t \\
&- \int_{\Omega}{ \Big( \partial_t \chi + \chi \Div_x \uvect \Big) B'(r) } 
\\ & + \int_{0}^{\tau} \int_{\partial \Omega} B'(r) ((r + \chi) - \rho_B) [\uvect_B \cdot \mathbf{n}]^- \ {\rm d} S_x \ \dt
\end{split}
\]
holds for any $0 \leq \tau \leq T$.

\end{lemma}

Using Lemma \ref{EL3}, we obtain strict positivity of $\rho$ on condition that $\rho_B$, $\rho_0$ have the same property. In fact, we have the following result.

\begin{coro}[{Abbatiello et al. \cite[Corollary 3.4]{Dissipative-Feireisl}}]\label{EC1}
Under the hypotheses of Lemma \ref{EL1}, we have 
\[
{\rm ess} \inf_{t,x} \rho(t,x) \geq  
\min \left\{ \min_\Omega \rho_0 , \,  \min_{\Gamma_{\rm in}} \rho_B \right\} 
\exp \left( -T \| \Div_x \uvect \|_{L^\infty((0,T) \times \Omega)} \right).
\]
\end{coro}

\subsection{Galerkin approximation of the momentum balance}   
Let us first define 
\begin{align}
\mathcal N_M = \left\{ \vc{v} : \vc{v} \in \C([0,T]; \C^2_c(\Omega; \mathbb R^3)), \ \|\vc{v}\|_{\C([0,T]; \C^2_c(\Omega; \mathbb R^3))} \leq M    \right\} ,
\end{align}
for some suitable constant $M$.
Then, we write the following result. 
\begin{lemma}\label{solution-Q}
For each $\vc{v} \in \C([0, T]; \C^{2}_c(\Omega, \R^3 ))$ with $\uvect=\vc{v}+\uvect_B$,
there exists a unique solution of the following initial-boundary value problem{\rm :}
\begin{align} \label{Q-system-FG}
\begin{cases}
\partial_{t} c + (\uvect\cdot \nabla_x ) c = D_{0} \Delta_x  c , \\
\partial_{t} Q + (\uvect \cdot \nabla_x ) Q + Q\Lambda-\Lambda Q 
= \Gamma H[Q, c] , \\
Q|_{t=0} = Q_{0} \in H^{1} (\Omega ) , \ \   Q |_{\partial  \Omega   } = Q_B \in H^{\sfrac{3}{2}}(\partial \Omega) , \\
c|_{t=0}=c_{0} \in H^{1}(\Omega), \quad   \nabla_x c\cdot  \mathbf n|_{\partial \Omega} = 0 ,
\end{cases}
\end{align} 
with $0<\underline{c} \leq c_{0}\leq \bar{c}<\infty$
such that
$Q\in L^{\infty}(0, T; H^{1}(\Omega))\cap L^{2}(0, T; H^{2}(\Omega))$, 
$c\in L^{\infty}(0, T ; H^{1}(\Omega))\cap L^{2}(0, T; H^{2}(\Omega))$,
and $0<\underline{c}\leq c \leq \bar{c}<\infty$. Moreover, the above mapping $\uvect\mapsto (c[\uvect], Q[\uvect])$ is continuous
from $\mathcal{N}_{M}$ to $\big(L^{\infty}(0, T; H^{1}(\Omega))\cap L^{2}(0, T; H^{2}(\Omega))\big)
 \times \big(L^{\infty}(0, T; H^{1}(\Omega))\cap L^{2}(0, T; H^{2}(\Omega))\big)$,
 and $Q[\uvect]\in S^{3}_{0}$ {\it a.e.} in $(0,T)\times \Omega$.  
\end{lemma} 
The proof of the above lemma is standard; one can follow a procedure similar to the one detailed in  \cite[Lemma 3.2]{C-Majumdar-W-Z-Comp}.



\vspace*{.2cm}

Now, the idea is to look for an approximate velocity field 
\[
\uvect = \vc{v} + \uvect_B , \ \ \vc{v} \in \C([0,T]; X_n) .
\] 
Then, the approximate momentum balance equation reads as
\begin{align}\label{approx-weak-momentum}
\int_\Omega \rho \uvect (\tau) \cdot \boldvphi(\tau) \d x &=
\int_0^\tau \int_\Omega \left[\rho \uvect \cdot \partial_t \boldvphi + \rho \uvect \otimes \uvect : \nabla_x \boldvphi + p(\rho) \Div_x \boldvphi  - \partial F_\delta(\D_x \uvect) : \nabla_x \boldvphi \right] \d x \dt \notag \\
& \ \  - \int_0^\tau \int_\Omega \left[G(Q)\mathbb I_3 - \nabla_x Q \odot \nabla_x Q \right] : \nabla_x \boldvphi \d x \d t  \notag \\
& \ \ - \int_0^\tau \int_\Omega \left[ Q \Delta_x Q - \Delta_x Q Q + \sigma_* c^2 Q   \right] : \nabla_x \boldvphi \d x \d t \notag \\ 
& \ \ - \veps \int_0^\tau \int_\Omega \nabla_x \rho \cdot \nabla_x \uvect \cdot \boldvphi \d x \d t 
+ \int_\Omega \rho_0 \uvect_0 \cdot \vphi(0) \dx ,
\end{align} 
for any $\boldvphi \in \C^1([0,T]; X_n)$, with the initial condition
\begin{equation} \label{E6}
(\rho \uvect)(0, \cdot)  =\rho_0 \uvect_0, \ \uvect_0 = \vc{v}_0 + \uvect_B, \ \vc{v}_0 \in X_n.
\end{equation}

For fixed parameter $n$, $\delta > 0$, $\veps > 0$, the first level approximation are solutions $(\rho, \uvect, c, Q)$ of the problem 
\eqref{cont-eps}--\eqref{initial-eps}, \eqref{Q-system-FG}, and the Galerkin approximation \eqref{approx-weak-momentum}. 
The \emph{existence} of the approximate solutions at this level can be proved in the same way as in \cite{ChJiNo}, see also \cite{C-Majumdar-W-Z-Comp}. Specifically, 
for given $\uvect =  \vc{v} + \uvect_B$, $\vc{v} \in \C([0,T]; X_n)$, we identify the unique solution $\rho = \rho[\uvect]$,  $c=c[\uvect]$, $Q=Q[\uvect]$ from the equations \eqref{cont-eps}--\eqref{initial-eps}, \eqref{Q-system-FG}. Then we plug them  into \eqref{approx-weak-momentum} which gives the unique solution $\uvect=\uvect[\rho, c, Q]$. To this end, we define a map
\begin{align*}
\mathcal T : \vc{v} \in \C([0,T]; X_n) \mapsto  \uvect [\rho, c, Q] - \uvect_B \in \C([0,T]; X_n) , 
\end{align*}
and the first level approximate solutions $\rho = \rho_{\delta, \veps, n}$, $\uvect = \uvect_{\delta, \veps, n}$, $c=c_{\delta, \veps, n}$ and $Q=Q_{\delta, \veps, n}$ are obtained via a fixed point argument, through the mapping $\mathcal{T}$. For more details about the process, we refer \cite{ChJiNo} (see also \cite{Sarka-Jan-Justyna}).

\vspace*{.2cm}

Finally, the approximate energy inequality is  
\begin{align}\label{Energy-ineq-approximate-delta}
&\frac{\d}{\d t} \int_{\Omega} \left(\frac{1}{2}\rho |\uvect - \uvect_B|^2 + P(\rho) + \frac{1}{2} |c|^2 + \frac{1}{2}|Q|^2 + \frac{1}{2}|\nabla_x Q|^2 + \frac{c_*}{4} |Q|^4 \right)   \d x  \notag \\
& + \frac{1}{2}\int_{\Omega} \partial F_\delta(\D_x\uvect) : \D_x \uvect \d x  
+  \frac{D_0}{2} \int_\Omega |\nabla_x c|^2 \d x \notag + \frac{\Gamma}{4} \int_{\Omega} |\Delta_x Q|^2  \d x  + 
\frac{c_*^2 \Gamma}{2} \int_{\Omega} |Q|^6\dx    \notag \\
& + \int_{\Gamma_{\out}} P(\rho) \uvect_B \cdot \mathbf n  \d S_x  + \int_{\Gamma_{\inn}} P(\rho_B) \uvect_B \cdot \mathbf n \d S_x  
\notag \\
&  -  \int_{\Gamma_{\inn}} \left[ P(\rho_B) - P'(\rho) (\rho_B - \rho) - P(\rho) \right] \uvect_B \cdot \vc{n} 
\d S_x   + \veps \int_\Omega P''(\rho) |\nabla_x\rho|^2 \d x   \notag \\
  \leq  & C \int_{\Omega} \left( \frac{1}{2} |c|^2 + \frac{1}{2}|Q|^2 + \frac{1}{2}|\nabla_x Q|^2 + \frac{c_*}{4} |Q|^4 \right)   \d x   + \int_\Omega \partial F_\delta(\D_x \uvect) : \nabla_x \uvect_B \dx \notag \\ 
  & - \int_{\Omega} \left[ \rho \uvect \otimes \uvect + p(\rho) \mathbb{I}_3\right] : \nabla_x \uvect_B \d x + \int_\Omega \rho \uvect \cdot (\uvect_B \cdot \nabla_x \uvect_B) \d x 
  \notag \\
 & - \frac{1}{2} \int_{\partial \Omega} \left(\frac{1}{2} |Q_B|^2 +  \frac{c_*}{4} \tr^2(Q_B^2) \right) \uvect_B \cdot \mathbf n   \d S_x 
  + 2c_* \Gamma \int_{\partial \Omega} |Q_B|^2 Q_B : \nabla_x Q_B \cdot \mathbf n \d S_x + C ,
\end{align}
where the constant $C>0$ may depend on $\uvect_B$, $Q_B$, $\Gamma$, and $D_0$  but independent of $\delta$, $\veps$ and $n$. 

\vspace{.2cm}
We conclude the above energy inequality by estimating the following quantities on the 
right-hand side  of \eqref{Energy-ineq-approximate-delta}. First, we compute 
\begin{align}\label{esti-1-p(rho)}
\left|\int_{\Omega} p(\rho) \mathbb{I}_3 : \nabla_x \uvect_B \d x \right| \leq C (\gamma-1)  \|\nabla_x \uvect_B\|_{L^\infty(\Omega)} P(\rho) ,
\end{align}
\begin{align}
\int_{\Omega} \rho \uvect \cdot (\uvect_B \cdot \nabla_x \uvect_B) \d x \leq C(\rho_0, \uvect_B) \left(1 + \int_{\Omega} \rho |\uvect- \uvect_B|^2 \d x \right) ,
\end{align}
and
\begin{align}  
\left|\int_{\Omega} \rho  \uvect \otimes \uvect :  \nabla_x \uvect_B \d x \right| \leq C  \|\nabla_x \uvect_B\|_{L^\infty(\Omega)} \left( \int_\Omega  \rho |\uvect - \uvect_B|^2 + C(\rho_0, \uvect_B) \right).
\end{align}
We also find that  
\begin{align}\label{esti-del-F-u}
\left|\int_{\Omega} \partial F_\delta (\mathbb D_x \uvect):\nabla_x \uvect_B \dx \right|&\leq
C \int_{\Omega} \sqrt{\partial F_\delta (\mathbb D_x \uvect):\mathbb D_x\uvect} \, |\nabla_x \uvect_B| \dx \notag \\
& \leq \frac{1}{4} \int_\Omega \partial F_\delta (\mathbb D_x \uvect):\mathbb D_x\uvect  \d x + C\|\nabla_x \uvect_B\|^2_{L^2(\Omega)}.
\end{align}

\vspace*{.1cm}

Then, after combining the estimates \eqref{esti-1-p(rho)}--\eqref{esti-del-F-u} in \eqref{Energy-ineq-approximate-delta}, we obtain 
\begin{align}\label{Energy-approx-del-final}
&\frac{\d}{\d t} \int_{\Omega} \left(\frac{1}{2}\rho |\uvect - \uvect_B|^2 + P(\rho) + \frac{1}{2} |c|^2 + \frac{1}{2}|Q|^2 + \frac{1}{2}|\nabla_x Q|^2 + \frac{c_*}{4} |Q|^4 \right)   \d x  \notag \\
& + \frac{1}{4}\int_{\Omega} \partial F_\delta(\D_x\uvect) : \D_x \uvect \d x  
+  \frac{D_0}{2} \int_\Omega |\nabla_x c|^2 \d x \notag + \frac{\Gamma}{4} \int_{\Omega} |\Delta_x Q|^2  \d x  + 
\frac{c_*^2 \Gamma}{2} \int_{\Omega} |Q|^6\dx    \notag \\
& + \int_{\Gamma_{\out}} P(\rho) \uvect_B \cdot \mathbf n  \d S_x   + \veps \int_\Omega P''(\rho) |\nabla_x\rho|^2 \d x  
\notag \\
&  -  \int_{\Gamma_{\inn}} \left[ P(\rho_B) - P'(\rho) (\rho_B - \rho) - P(\rho) \right] \uvect_B \cdot \vc{n} 
\d S_x     \notag \\
  \leq  & C \int_{\Omega} \left(\frac{1}{2}\rho |\uvect - \uvect_B|^2 + P(\rho) + \frac{1}{2} |c|^2 + \frac{1}{2}|Q|^2 + \frac{1}{2}|\nabla_x Q|^2 + \frac{c_*}{4} |Q|^4 \right)   \d x \notag \\ 
 & - \frac{1}{2} \int_{\partial \Omega} \left(\frac{1}{2} |Q_B|^2 +  \frac{c_*}{4} \tr^2(Q_B^2) \right) \uvect_B \cdot \mathbf n   \d S_x  -  \int_{\Gamma_{\inn}} P(\rho_B) \uvect_B \cdot \mathbf n \d S_x \notag \\ 
&  +  2c_* \Gamma \int_{\partial \Omega} |Q_B|^2 Q_B : \nabla_x Q_B \cdot \mathbf n \d S_x + C ,
\end{align}
where the constant $C>0$  is independent of $\delta$, $\veps$ and $n$.

\begin{remark}
   Note that, in the left-hand side of the  energy inequality \eqref{Energy-approx-del-final},  the term $\disp \int_{\Gamma_{\out}}P(\rho) \uvect_B \cdot \mathbf{n} \d S_x$ has positive sign since $\uvect_B \cdot \mathbf{n} \geq 0$ in $\Gamma_{\out}$.  

   On the other hand, since $\uvect_B\cdot \mathbf{n}<0$ in $\Gamma_{\inn}$ and due to the convexity of $P(\rho)$, we have  $\disp \left[ P(\rho_B) - P'(\rho) (\rho_B - \rho) - P(\rho) \right]\geq 0$,  the associated term in the left-hand side of \eqref{Energy-approx-del-final} is non-negative. These facts guarantee the validity of the approximate energy inequality \eqref{Energy-approx-del-final}.  
\end{remark}

\vspace*{.2cm}

This leads to the following proposition. 
\begin{prop}[\textbf{Solution to approximate system}]\label{Proposition-Approx-I}
Let $\Omega \subset \R^3$ be a bounded Lipschitz domain. Suppose that $p = p(\rho)$ as in \eqref{S2} and the regularized potential  $F_\delta$ satisfies \eqref{F-delta}--\eqref{F-delta-convex}. Let the data belong to the class 
\begin{align*} 
\uvect_B \in \C^{1}_c(\R^3; \R^3), \ \rho_{0} \in \C^{1}(\R^3) \text{ such that } \rho_0\equiv 0 \text{ in } \mathbb R^3\setminus \Omega , \   \rho_B \in \C^1(\partial \Omega),\ \rho_0, \rho_B \geq   0, \\
\uvect_0 = \mathbf v_0 + \uvect_B,\ \mathbf v_0 \in X_n, \ c_0 \in \C^1_c(\Omega), \ Q_0 \in \C^1_c(\Omega), \ Q_B \in \C^1(\partial \Omega) .
\end{align*}

Then for each fixed $\delta > 0$, $\veps > 0$ and $n > 0$, there exists a solution 
$(\rho, \uvect, c, Q)$ of the approximate problem \eqref{cont-eps}--\eqref{initial-eps}, \eqref{Q-system-FG} and \eqref{approx-weak-momentum}. Moreover, the approximate energy inequality \eqref{Energy-approx-del-final} holds. 
\end{prop}


\section{Limit with respect to $\delta \to 0$}\label{Section-delta-to-0}

In this section, we pass to the $\delta\to 0$ limit in the regularization of the  potential $F_\delta$. Firstly,
we apply 
 Gr\"{o}nwall's lemma in \eqref{Energy-approx-del-final}, and as a consequence we have the following estimates: 
\begin{align}\label{Energy-ineq-Gronwall}
& \left[\int_{\Omega} \left(\frac{1}{2}\rho |\uvect - \uvect_B|^2 + P(\rho) +  \frac{1}{2}|c|^2  + \frac{1}{2}|Q|^2  + \frac{1}{2}|\nabla_x Q|^2  + \frac{c_*}{4}|Q|^4 \right)   \d x\right]_{0}^{\tau}  \notag \\
& 
+  \int_0^\tau \left(\frac{D_0}{2}\|\nabla 
_x c\|^2_{L^2(\Omega)}  + \frac{\Gamma}{4} \|\Delta_xQ\|^2_{L^2(\Omega)} + \frac{c_*^2 \Gamma}{2}\|Q\|^{6}_{L^6(\Omega)} \right) \d t   \notag \\
& + \frac{1}{4}\int_0^\tau \int_{\Omega} \partial F_\delta (\D_x \uvect): \D_x \uvect \d x \d t    + \int_0^\tau \int_{\Gamma_{\out}} P(\rho) \uvect_B \cdot \mathbf n  \d S_x \d t 
\notag  \\
& - \int_0^\tau \int_{\Gamma_{\inn}} \left[ P(\rho_B) - P'(\rho) (\rho_B - \rho) - P(\rho) \right] \uvect_B \cdot \vc{n} 
\d S_x \d t  + \veps \int_0^\tau \int_\Omega P''(\rho) |\nabla_x\rho|^2 \d x  \d t 
\notag \\
 &\leq  C , 
\end{align}
for a.a. $\tau \in [0,T]$, where $C>0$ is a constant that depends on given data and parameters  but independent of $\delta$, $\veps$ and $n$. 

\vspace*{.2cm} 

Thanks to \eqref{Energy-ineq-Gronwall} and the relation \eqref{F-delta-convex},  we deduce the uniform bound on the traceless part of the velocity gradient. Indeed we have  
\begin{align*}
    \left\| \D_x \uvect_\delta - \frac{1}{3}\Div_x \uvect_\delta \mathds{I}_3  \right\|^{\sfrac{4}{3}}_{L^{\sfrac{4}{3}}(0,T; L^{\sfrac{4}{3}}(\Omega; \R^{3\times 3}) )} \lesssim \int_0^T \int_\Omega \partial F_\delta (\D_x \uvect_\delta): \D_x \uvect_\delta \d x \d t \leq C ,
\end{align*}
uniformly w.r.t. $\delta \to 0$.
Then, one may use the trace-free Korn's inequality given by Lemma \ref{Lemma-Korn-trace-free},  to deduce that    
\begin{align*}
\|\nabla_x \uvect_\delta \|_{L^{\sfrac{4}{3}}(0,T; L^{ \sfrac{4}{3} }(\Omega; \R^{3\times 3}))}\leq C . 
\end{align*}
Combining with the standard Poincar\'{e}  inequality,  we further obtain that 
\begin{align}\label{bound-u-delta}
\|\uvect_\delta \|_{L^{\sfrac{4}{3}}(0,T; W^{1,\sfrac{4}{3}}(\Omega; \R^3) )} \leq C .
\end{align}
Since we are in finite dimensional space $X_n$, all  norms are equivalent, in turn we have 
\begin{align}\label{bound-u-del-L-infty}
\|\uvect_\delta\|_{L^\infty(0,T; X_n)} + \|\nabla_x \uvect_\delta\|_{L^\infty(0,T; X_n)} \leq \|\uvect_\delta \|_{L^{\sfrac{4}{3}}(0,T; W^{1,\sfrac{4}{3}}(\Omega; \R^3))} \leq C . 
\end{align}

\vspace*{.2cm}

For the density $\rho_\delta$, we have, in  view of Lemma \ref{EL2} and Corollary \ref{EC1}, the following uniform bounds: 
\begin{equation} \label{bound-rho-up-low}
0 < \underline{\rho} \leq \rho_\delta(t,x) \leq \overline{\rho} \  \ \  \forall  (t,x) \in [0,T] \times \overline{\Omega}. 
\end{equation} 
Moreover, from \eqref{Energy-ineq-Gronwall}, we  obtain 
\begin{align}\label{bound-rho-u-del}
\|\rho_\delta |\uvect_\delta|^2\|_{L^\infty(0,T ; L^{1}(\Omega))} + \|\rho_\delta \|_{L^\infty(0,T; L^{\gamma}(\Omega))} \leq C .
\end{align}
At this stage, since the approximated continuity equation has parabolic effect, we can gain more information about the density $\rho_\delta$. More precisely, starting with regular enough initial density $\rho_0$, one can obtain the following bounds (taking into account the uniform bound \eqref{bound-u-del-L-infty})
\begin{align}
\|\rho_\delta\|_{L^2(0,T; H^2(\Omega))} + \|\partial_t\rho_\delta\|_{L^2((0,T)\times\Omega)}\leq C,
\end{align}
which yields  (up to a subsequence, if necessary) 
\begin{align}\label{limit-rho-delta-1}
\begin{dcases} 
\rho_\delta \to \rho \text{ weakly in } L^2(0,T; H^2(\Omega)) , \\ 
\partial_t \rho_\delta \to \rho \text{ weakly in } L^2((0,T)\times \Omega).
\end{dcases} 
\end{align} 
Consequently,  by Aubin-Lions lemma, 
we get
\begin{align}\label{limit-strong-rho-delta}
\rho_\delta \to \rho \text{ strongly in } L^2(0,T; H^1(\Omega)) . 
\end{align}
The above  result \eqref{limit-strong-rho-delta} further yields the convergences of  the pressure term and its potential, namely we have 
\begin{align}\label{conv-pressure-delta}
p(\rho_\delta) \to p(\rho) \ \text{ and } \ P(\rho_\delta) \to P(\rho) \ \text{ strongly in } L^r((0,T)\times \Omega) 
\end{align}
for some $r>1$. 

\vspace*{.2cm}

Next, from the energy inequality \eqref{Energy-ineq-Gronwall}, one can also conclude that
\begin{align}\label{bound-c-Q-del}
&\| c_\delta \|_{L^\infty(0,T; L^{2}(\Omega))} + 
\| c_\delta \|_{L^2(0,T; H^1(\Omega))}  \notag  \\
&+ \| Q_\delta\|_{L^\infty(0,T; H^1(\Omega))} + 
\| Q_\delta\|_{L^2(0,T; H^2(\Omega))} + \|Q_\delta\|_{L^6((0,T)\times \Omega)}\leq C  .
\end{align}
In turn, we have 
\begin{align}\label{conv-c-Q-delta}
c_\delta \to c \text{ weakly in } L^2(0,T; H^1(\Omega)) \ \text{ and } \ Q_\delta \to Q \text{ weakly in } L^2(0,T; H^2(\Omega))  . 
\end{align}
One can further find the uniform bounds of $\partial_t c_\delta$ and $\partial_t Q_\delta$ at least in the spaces $L^1(0,T; H^{-1}(\Omega))$ and $L^r(0,T; L^q(\Omega))$ for some $r, q>1$ (by performing a similar analysis which we will describe later in Section \ref{Section-Epsilon}). As a consequence, by using Aubin-Lions lemma, we find 
\begin{align}\label{strong-conv-c-Q-del}
c_\delta \to c \text{ strongly in } L^2((0,T)\times \Omega) \ \text{ and } \ Q_\delta \to Q \text{ strongly in } L^2(0,T; H^1(\Omega)) . 
\end{align}

\vspace{.2cm}

Once we have all the uniform bounds concerning $\rho_\delta, c_\delta, Q_\delta$,  one can follow the steps in  \cite[Section 5]{Sarka-Jan-Justyna} to bound the time derivative of $\uvect_\delta$ as follows:   
\begin{align*}
\|\partial_t \uvect_\delta\|_{L^2(0,T; X_n)} \leq C ,
\end{align*}
where $C>0$ is uniform in $\delta$. 
With this bound in hand, along with \eqref{bound-u-del-L-infty}, yields the following strong convergence result:  
\begin{align}\label{strong-conv-u-delta}
\uvect_\delta \to  \uvect \text{ strongly in } \C([0,T]; X_n) , 
\end{align}
 thanks to Aubin-Lions lemma. 

\vspace{.2cm} 

Finally, we need to establish   
the convergence of the term $\partial F_\delta(\D_x\uvect_\delta)$ in the approximate momentum equation \eqref{approx-weak-momentum}, which we shall analyze below.

\vspace*{.2cm}

\noindent 
-- {\em Limit of $\partial F_\delta(\D_x \uvect_\delta)$.}
First, note that
\begin{align*}
F_\delta(\D) \leq G(\D): = \sup_{\delta \in [0,1]} F_\delta (\D) .
\end{align*}
Since the functions $F_\delta$ are proper convex, the same holds for $G$ on $\mathbb R^{3\times 3}_{\sym}$ and moreover, its conjugate $G^*$ is superlinear, and 
\begin{align}\label{inequality-F-G-star}
G^*(\D) \leq F_{\delta}^*(\D) \quad \forall \D \in \mathbb R^{3\times 3}_{\sym} , \ \delta \in [0,1] .
\end{align}
Then, using the relation \eqref{stress-tensor}--\eqref{relation-stress-tensor} between subgradients and conjugate functions of the smooth function $F_\delta$, we write 
\begin{align*}
\int_0^T \int_{\Omega} \partial F_\delta (\D_x \uvect_\delta) : \D_x \uvect_\delta   \d x \d t = 
\int_0^T \int_{\Omega} \left[F_\delta (\D_x \uvect_\delta) + F^* (\partial F_\delta (\D_x \uvect_\delta) )\right]  \d x \d t .
\end{align*}
Using the energy estimate \eqref{Energy-ineq-Gronwall} and \eqref{inequality-F-G-star}, we deduce 
\begin{align*}
\int_0^T \int_{\Omega} G^*(\partial F_\delta(\D_x \uvect_\delta)) \d x \d t &\leq \int_0^T \int_{\Omega} F^*_{\delta}(\partial F_\delta(\D_x \uvect_\delta)) \d x \d t \\
&\leq \int_0^T \int_{\Omega} \partial F_\delta (\D_x \uvect_\delta) : \D_x \uvect_\delta   \d x \d t\leq C ,
\end{align*}
for some constant $C>0$ which is independent in $\delta$. 

 Then, due to the superlinearity of $G^*$, one may apply De La Vall\'{e}e--Poussin criterion for equi-integrability (see Lemma 
 \ref{Lemma-La-Valle}) to deduce that there exists a function $\mathbb{S}\in L^1((0,T)\times \Omega; \mathbb R^{3\times 3}_{\sym})$ such that
\begin{align}\label{weak-conv-subgradient}
\partial F_\delta(\D_x \uvect_\delta) \to  \mathbb{S} \, \text{ weakly in } L^1((0,T)\times \Omega).
\end{align}
One now needs to identify $\mathbb{S}$ with the subgradient of $F$ at $\D_x \uvect$ (i.e., $\partial F(\D_x \uvect)$), where $\uvect$ is the limiting function given in \eqref{strong-conv-u-delta}.  
We achieve this goal by following
the arguments given in \cite[Section 3.1]{Woźnicki2021mv}. Note that,  $\partial F_\delta(\D_x \uvect_\delta)$  constitutes the classical gradient (as it is smooth by construction \eqref{F-delta}), which is 
in particular a subgradient of $F_\delta$ at $\D_x \uvect_\delta$. Thus, for all 
non-negative functions $\phi\in \mathcal D((0,T)\times \Omega)$ and all matrices $\D\in \R^{3\times 3}_{\sym}$,  it holds that
\begin{align}\label{relation-subgradient}
\int_0^T \int_\Omega \Big[F_\delta(\D)  - \partial F_{\delta}(\D_x \uvect_\delta) : (\D - \D_x \uvect_\delta) - F_\delta(\D_x \uvect_\delta)\Big] \phi \d x \d t \geq 0 . 
\end{align}
Further, from the definition \eqref{F-delta} of the regularized potential $F_\delta$, we know that
\begin{align}
F_\delta(\cdot) \to F(\cdot) \ \text{ in } \C_{\textnormal{loc}}(\R^{3\times 3}_{\sym}) .
\end{align}
This, combining  with the uniform convergence of $\D_x \uvect_\delta$ using \eqref{bound-u-del-L-infty} and the weak convergence \eqref{weak-conv-subgradient} allows
us to pass to the limit on the right-hand side of the inequality \eqref{relation-subgradient}. In fact, due to the arbitrary choice of $\phi\in \mathcal D((0,T)\times \Omega)$, one has
\begin{align}
    F(\D)  - \Ss : (\D - \D_x \uvect) - F(\D_x \uvect) \geq 0 \quad \forall \, \D \in \R^{3\times 3}_{\sym}  ,
\end{align}
which is the desired criterion for $\Ss\in \partial F(\D_x \uvect)$. In particular, the following relations hold:
\begin{align}
\Ss\in \partial F(\D_x\uvect) \Leftrightarrow \D_x \uvect \in \partial F^*(\Ss) \Leftrightarrow 
\Ss : \D_x \uvect = F(\D_x\uvect) + F^*(\Ss). 
\end{align}

\vspace*{.15cm}

With   the above convergence results, we can now pass  to the limit  w.r.t. $\delta\to 0$ to obtain the weak formulations satisfied by the limiting functions  $(\rho, \uvect, c, Q)$ and $\Ss$.  

\subsection{Limits in the continuity, concentration and tensor parameter equations}

After passing to the $\delta\to 0$ limit, the continuity equation, the equation for active particles, and the equation for  nematic tensor order parameter  have similar formulations as in \eqref{E4}, \eqref{weak-concen} and \eqref{weak-nematic-tensor} respectively, thanks to the convergence results obtained in this level. We do not repeat those equations. 


\subsection{Limit in the momentum balance}
In accordance with all the limits obtained above w.r.t. $\delta\to 0$, the momentum 
equation \eqref{approx-weak-momentum}  now reads as
\begin{align}\label{weak-momentum-delta-limit}
\int_\Omega \rho \uvect (\tau) \cdot \boldvphi(\tau) \d x &=
\int_0^\tau \int_\Omega \left[\rho \uvect \cdot \partial_t \boldvphi + \rho \uvect \otimes \uvect : \nabla_x \boldvphi + p(\rho) \Div_x \boldvphi  - \Ss : \nabla_x \boldvphi \right] \d x \dt \notag \\
& - \int_0^\tau \int_\Omega \left[G(Q)\mathbb I_3 - \nabla_x Q \odot \nabla_x Q \right] : \nabla_x \boldvphi \d x \d t  \notag \\
&- \int_0^\tau \int_\Omega \left[ Q \Delta_x Q - \Delta_x Q Q + \sigma_* c^2 Q   \right] : \nabla_x \boldvphi \d x \d t \notag \\ 
&- \veps \int_0^\tau \int_\Omega \nabla_x \rho \cdot \nabla_x \uvect \cdot \boldvphi \d x \d t 
+ \int_\Omega \rho_0 \uvect_0 \cdot \vphi(0) \dx ,
\end{align} 
for any $\boldvphi \in \C^1([0,T]; X_n)$, with the initial condition
\begin{equation*}
(\rho \uvect)(0, \cdot)  =\rho_0 \uvect_0, \ \uvect_0 = \vc{v}_0 + \uvect_B, \ \vc{v}_0 \in X_n.
\end{equation*}

\subsection{Limiting energy inequality}
Finally, the energy inequality for the limiting system (as $\delta \to 0$) is
\begin{align}\label{Main-Energy-limit-delta}
&\frac{\d}{\d t} \int_{\Omega} \left(\frac{1}{2}\rho |\uvect - \uvect_B|^2 + P(\rho) + \frac{1}{2} |c|^2 + \frac{1}{2}|Q|^2 + \frac{1}{2}|\nabla_x Q|^2 + \frac{c_*}{4} |Q|^4 \right)   \d x  \notag \\
& + \frac{1}{2}\int_{\Omega} \left[ F(\D_x\uvect) + F^*(\Ss) \right] \d x  
+  \frac{D_0}{2} \int_\Omega |\nabla_x c|^2 \d x \notag + \frac{\Gamma}{4} \int_{\Omega} |\Delta_x Q|^2  \d x  + 
\frac{c_*^2 \Gamma}{2} \int_{\Omega} |Q|^6\dx    \notag \\
& + \int_{\Gamma_{\out}} P(\rho) \uvect_B \cdot \mathbf n  \d S_x  + \int_{\Gamma_{\inn}} P(\rho_B) \uvect_B \cdot \mathbf n \d S_x  
\notag \\
&  -  \int_{\Gamma_{\inn}} \left[ P(\rho_B) - P'(\rho) (\rho_B - \rho) - P(\rho) \right] \uvect_B \cdot \vc{n} 
\d S_x   + \veps \int_\Omega P''(\rho) |\nabla_x\rho|^2 \d x   \notag \\
  \leq  & C \int_{\Omega} \left( \frac{1}{2} |c|^2 + \frac{1}{2}|Q|^2 + \frac{1}{2}|\nabla_x Q|^2 + \frac{c_*}{4} |Q|^4 \right)   \d x   + \int_\Omega \Ss : \nabla_x \uvect_B \dx   \notag \\ 
  & - \int_{\Omega} \left[ \rho \uvect \otimes \uvect + p(\rho) \mathbb{I}_3\right] : \nabla_x \uvect_B \d x + \int_{\Omega} \rho \uvect \cdot (\uvect_B \cdot \nabla_x \uvect_B) \d x 
 \notag \\
 & - \frac{1}{2} \int_{\partial \Omega} \left(\frac{1}{2} |Q_B|^2 +  \frac{c_*}{4} \tr^2(Q_B^2) \right) \uvect_B \cdot \mathbf n   \d S_x 
 + 2c_* \Gamma \int_{\partial \Omega} |Q_B|^2 Q_B : \nabla_x Q_B \cdot \mathbf n \d S_x+ C ,
\end{align}
where the constant is independent  of $\veps$ or $n$.

Accordingly, the energy estimate is  (after applying the Gr\"{o}nwall's lemma as earlier) 
\begin{align}\label{Energy-ineq-delta-limit}
& \left[\int_{\Omega} \left(\frac{1}{2}\rho |\uvect - \uvect_B|^2 + P(\rho) +  \frac{1}{2}|c|^2  + \frac{1}{2}|Q|^2  + \frac{1}{2}|\nabla_x Q|^2  + \frac{c_*}{4}|Q|^4 \right)   \d x\right]_{0}^{\tau}  \notag \\
& 
+  \int_0^\tau \left(\frac{D_0}{2}\|\nabla 
_x c\|^2_{L^2(\Omega)}  + \frac{\Gamma}{4} \|\Delta_xQ\|^2_{L^2(\Omega)} + \frac{c_*^2 \Gamma}{2}\|Q\|^{6}_{L^6(\Omega)} \right) \d t   \notag \\
& + \frac{1}{4}\int_0^\tau \int_{\Omega} \left[ F (\D_x \uvect) + F^*(\Ss) \right] \d x \d t    + \int_0^\tau \int_{\Gamma_{\out}} P(\rho) \uvect_B \cdot \mathbf n  \d S_x \d t 
\notag  \\
& - \int_0^\tau \int_{\Gamma_{\inn}} \left[ P(\rho_B) - P'(\rho) (\rho_B - \rho) - P(\rho) \right] \uvect_B \cdot \vc{n} 
\d S_x \d t  + \veps \int_0^\tau \int_\Omega P''(\rho) |\nabla_x\rho|^2 \d x  \d t 
\notag \\
 &\leq  C 
\end{align}
for a.a. $\tau \in [0,T]$, where the constant $C>0$ is independent of $\veps$ and $n$.

\vspace*{.2cm}

Combining the estimates and results above, we state the proposition below. 
\begin{prop}\label{Prop-Approx-II}
Assume that the hypothesis in Proposition \ref{Proposition-Approx-I} hold true. Let $(\rho_\delta, \uvect_\delta, c_\delta, Q_\delta)_{\delta>0}$ be the sequence of solutions to the approximate problem \eqref{cont-eps}--\eqref{initial-eps}, \eqref{Q-system-FG} and \eqref{approx-weak-momentum} for given $\veps>0$ and $n\in \mathbb N$. Then, there exist some functions $(\rho, \uvect, c, Q)$ such that the following convergence results hold:
\begin{align*}
& \rho_\delta \to \rho  \, \text{ strongly in } \, L^2(0,T; H^1(\Omega)) , \ \ \uvect_\delta \to \uvect \, \text{ strongly in } \, \C([0,T]; X_n)  , \\
&  c_\delta \to c \, \text{ weakly in } \, L^2(0,T; H^1(\Omega))  , \ \ c_\delta \to c \, \text{ strongly in } \, L^2((0,T)\times \Omega) , \\
&  Q_\delta \to Q \, \text{ weakly in } \, L^2(0,T; H^2(\Omega)) , \ \ Q_\delta \to Q \, \text{ strongly in } \, L^2(0,T; H^1(\Omega))  . 
\end{align*}
In addition,  there exists a function $\mathbb{S}\in L^1((0,T)\times \Omega; \mathbb R^{3\times 3}_{\sym})$ such that
\begin{align*}
\partial F_\delta(\D_x \uvect_\delta) \to  \mathbb{S} \, \text{ weakly in } \, L^1((0,T)\times \Omega), 
\end{align*}
where $\Ss$ satisfies 
\begin{align*}
\Ss : \D_x \uvect = F(\D_x\uvect) + F^*(\Ss) ,
\end{align*}
and $F$ is given by \eqref{Coercivity-F}. 

Moreover, $(\rho, \uvect, c, Q, \Ss)$ satisfy the weak formulations presented in this section.  
\end{prop}

\vspace*{.2cm}

 Thus, we end up with  
 the $\veps$-level solutions. In fact,  the sequence of solutions  obtained in Proposition \ref{Prop-Approx-II} may be re-denoted as $(\rho_\veps, \uvect_\veps, c_\veps, Q_\veps)_{\veps>0}$ and $(\Ss_\veps)_{\veps>0}$. We shall treat them  in the next section.



\section{Limit with respect to $\varepsilon \to 0$}\label{Section-Epsilon}

This section is devoted to the $\veps\to 0$ limit
 in the viscous approximation \eqref{cont-eps}.  
This is more intricate than the preceding step as we  lose compactness of the approximate density w.r.t.  the spatial variable. For fixed $n > 0$,  we collect the necessary estimates  that are independent of $\veps$. 

\vspace*{.15cm} 

Similar to Section \ref{Section-delta-to-0} (as we obtained \eqref{bound-u-delta}), we have 
\begin{equation}\label{bound-u-eps}
\| \uvect_\veps \|_{L^{\sfrac{4}{3}}(0,T; W^{1,\sfrac{4}{3}}(\Omega;\mathbb R^3) )} \leq C, 
\end{equation}
thanks to the energy estimate \eqref{Energy-ineq-delta-limit}. Since $n$ is still fixed, we have 
\begin{align}\label{bound-u-eps-finite}
\|\uvect_\veps\|_{L^{\sfrac{4}{3}}(0,T; X_n)} \leq \| \uvect_\veps \|_{L^{\sfrac{4}{3}}(0,T; W^{1,\sfrac{4}{3}}(\Omega;\mathbb R^3) )} \leq C.
\end{align}
This yields, in view of Lemma \ref{EL2} and Corollary \ref{EC1}, the following uniform bounds on the density $\rho_\veps$:
\begin{equation} \label{bound-rho-up-low}
0 < \underline{\rho} \leq \rho_\veps(t,x) \leq \overline{\rho} \  \ \  \forall  (t,x) \in [0,T] \times \overline{\Omega}.
\end{equation} 
The energy estimate \eqref{Energy-ineq-delta-limit} further yields 
\begin{align}\label{bound-sqrt-rho-u-eps}
\|\rho_\veps |\uvect_\veps|^2\|_{L^\infty(0,T; L^1(\Omega))} + \|\rho_\veps \|_{L^\infty(0,T; L^\gamma(\Omega))} \leq C . 
\end{align}
Moreover, from  \eqref{Energy-ineq-delta-limit}, the upper and lower bounds \eqref{bound-rho-up-low}, and \eqref{bound-u-eps-finite}, we infer that
\begin{align}\label{L-infty-u-bound}
\esssup_{\tau \in [0,T]} \|\uvect_\veps (\tau)\|_{W^{1,\infty}(\Omega;\R^3)} \leq C .
\end{align}

Note that at this stage $\rho_\veps$ possesses a 
well-defined trace on $\partial \Omega$. Thus, testing the equation \eqref{cont-eps} with $\rho_\veps$ and using the bounds \eqref{bound-u-eps}--\eqref{bound-rho-up-low}, one has 
\begin{equation} \label{bound-grad-rho-eps}
\veps \int_0^T \int_{\Omega}{|\nabla_x \rho_\veps|^2 } \dt \leq C, 
\end{equation}
where $C$ is independent in $\veps$.  

\vspace*{.15cm} 

In view of the uniform bounds obtained before, one may assume 
\begin{align}\label{limit-rho-eps-1}
\begin{dcases}
\rho_{\veps} \to \rho \ \mbox{weakly-(*) in}\ L^\infty((0,T) \times \Omega) \ \mbox{and} \\
\rho_{\veps} \to \rho   \ \mbox{weakly in}\ \C_{\rm weak}([0,T]; L^r(\Omega)) \ \mbox{for any}\ 1 < r < \infty, 
\end{dcases}
\end{align}
up to a suitable subsequence if necessary. The second convergence follows form the weak bound on the 
time derivative $\partial_t \rho_{\veps}$ obtained from equation \eqref{E4}. We also have 
\begin{align}\label{limit-boundary-rho-eps}
\rho_{\veps} \to \rho \ \mbox{weakly-(*) in}\ L^\infty((0,T) \times \partial \Omega ;  \d S_x).
\end{align} 
Additionally, the limit density admits the same upper and lower bounds as in \eqref{bound-rho-up-low}.

\vspace*{.2cm} 

Next, from the uniform estimate \eqref{L-infty-u-bound}, we have  
\begin{equation} \label{E15b}
\uvect_\veps \to \uvect \ \mbox{weakly-(*) in} \ L^\infty(0,T; W^{1,\infty}(\Omega; \R^3)), 
\end{equation}
and 
\begin{align}\label{limit-rhoU-eps}
\rho_{\veps} \uvect_{\veps} \to \mathbf{m} \ \mbox{weakly-(*) in} \ L^\infty((0,T) \times \Omega; \R^3)).
\end{align} 
Then, by using Arzela--Ascoli theorem, one can show that 
\begin{equation} \label{E15c}
\mathbf{m} = \rho \uvect \ \, \mbox{a.e. in}\ \, (0,T) \times \Omega.
\end{equation}

\vspace*{.2cm} 

Furthermore, from the uniform energy estimate \eqref{Energy-ineq-delta-limit} and by using the maximum  principle  to the equation of $c_\veps$, we deduce that  
\begin{align}\label{bound-c-Q-eps}
\begin{dcases} 
\|c_\veps\|_{L^\infty((0,T)\times \Omega)} + \|c_\veps\|_{L^2(0,T; H^1(\Omega))}  \leq C , \\
\|Q_\veps \|_{L^\infty(0,T; H^1(\Omega))}  + \|Q_\veps \|_{L^2(0,T;H^2(\Omega))} \leq C ,
\end{dcases}
\end{align} 
where $C>0$ is independent of $\veps$. Consequently, we have  
\begin{align}\label{conv-c-Q-eps}
c_\veps \to c \ \mbox{ weakly in } L^2(0,T; H^1(\Omega)) \ \mbox{ and } \ Q_\veps \to Q \ \mbox{ weakly in } L^2(0,T; H^2(\Omega)) .
\end{align}

\vspace*{.2cm} 

We  still need to find uniform bounds of the terms $\partial_t c_\veps$ and $\partial_t Q_\veps$. To deduce the former one, let us recall the equation \eqref{concen-eq-1}, one has 
\begin{align*}
\int_0^T \left\langle \partial_t c_\veps , \phi\right\rangle_{H^{-1}, H^1_0} \d t = -\int_0^T \left\langle (\uvect_\veps  \cdot \nabla_x) c_\veps  ,  \phi\right\rangle_{H^{-1}, H^1_0} \d t + D_0\int_0^T \left\langle \Delta_x c_\veps , \phi\right\rangle_{H^{-1}, H^1_0} \d t ,
\end{align*}
for all $\phi \in H^1_0(\Omega)$ with $\|\phi\|_{H^1_0(\Omega)} = 1 . $
Then, we compute 
\begin{align*}
&- \left\langle (\uvect_\veps \cdot \nabla_x) c_\veps  ,  \phi\right\rangle_{H^{-1}, H^1_0}  + D_0 \left\langle \Delta_x c_\veps , \phi\right\rangle_{H^{-1}, H^1_0}  \\
 =& - \int_\Omega (\uvect_\veps \cdot \nabla_x) c_\veps  \,  \phi\d x - D_0 \int_\Omega \nabla_x c_\veps \cdot \nabla_x \phi \d x  \\
 =&  \int_{\Omega} c_\veps \, \Div_x(\uvect_\veps \phi) \d x - D_0\int_\Omega \nabla_x c_\veps \cdot \nabla_x \phi \d x , 
\end{align*}
which yields 
\begin{align*}
\Big|\left\langle \partial_t c_\veps , \phi\right\rangle_{H^{-1}, H^1_0} \Big| 
& 
\leq \|c_\veps\|_{L^\infty(\Omega)} \left( \|\nabla_x \uvect_\veps\|_{L^{\sfrac{4}{3}}(\Omega; \R^{3\times 3})} \|\phi\|_{L^{4}(\Omega)}     
+ \|\uvect_\veps\|_{L^{\sfrac{12}{5}}(\Omega; \R^3)} \|\nabla_x \phi\|_{L^{\sfrac{12}{7}}(\Omega)} \right) 
\\ &  \quad + D_0 \|\nabla_x c_\veps \|_{L^2(\Omega)}  \|\nabla_x \phi \|_{L^2(\Omega)} \\
& \leq C \|c_\veps\|_{L^\infty(\Omega)} \|\uvect_\veps \|_{W^{1,\sfrac{4}{3}}(\Omega; \R^3)} \|\phi\|_{H^1_0(\Omega)} + D_0 \|c_\veps \|_{H^1(\Omega)}  \|\phi \|_{H^1_0(\Omega)}  ,  
\end{align*}
thanks to the compact embeddings $W^{1,\sfrac{4}{3}}(\Omega) \hookrightarrow L^{\sfrac{12}{5}}(\Omega)$ and $H^1_0(\Omega)\hookrightarrow L^4(\Omega)$ (recall that we are in spatial dimension $3$).

Therefore, we have 
\begin{align}\label{bound-c-t-eps}
\|\partial_t c_\veps \|_{L^1(0,T; H^{-1}(\Omega))} &\leq C \|c_\veps\|_{L^\infty((0,T)\times \Omega)} \|\uvect_\veps\|_{L^{\sfrac{4}{3}}(0,T; W^{1,\sfrac{4}{3}}(\Omega;\R^3)  } + D_0 \|c_\veps\|_{L^2(0,T;H^1(\Omega))}  \notag  \\
& \leq  C,
\end{align}
by using the uniform bounds \eqref{bound-u-eps} and \eqref{bound-c-Q-eps}.

\vspace*{.2cm}

We shall now find  proper estimate for the quantity  $\partial_t Q_\veps$. Recall the equation of $Q_\veps$ from \eqref{nematic-tensor-eq-1}, and we find the following bounds:
\begin{align}\label{Q-t-1}
\|(\uvect_\veps \cdot \nabla_x) Q_\veps\|_{L^{\sfrac{4}{3}}(0,T; L^{\sfrac{12}{11}}(\Omega))} 
 &\leq \bigg(\int_{0}^{T}\|\uvect_\veps\|^{\sfrac{4}{3}}_{L^{\sfrac{12}{5}}(\Omega;\R^3)} \|\nabla_x Q_\veps \|_{L^{2}(\Omega)}^{\sfrac{4}{3}} \d t \bigg)^{\sfrac{3}{4}}  \notag \\
 & \leq C \|\uvect_\veps\|_{L^{\sfrac{4}{3}}(0,T; W^{1,\sfrac{4}{3}}(\Omega;\R^3)) }\|\nabla_x Q_\veps\|_{L^{\infty}(0,T; L^{2}(\Omega))} 
 \notag \\ 
 & \leq C, 
 \end{align}
 \begin{align}\label{Q-t-2} 
\|Q_\veps\Lambda_\veps - \Lambda_\veps Q_\veps\|_{L^{\sfrac{4}{3}}(0,T;L^{\sfrac{12}{11}}(\Omega))} 
&\leq 
2 \bigg(\int_0^T \|\nabla_x \uvect_\veps\|^{\sfrac{4}{3}}_{L^{\sfrac{4}{3}}(\Omega; \R^{3\times 3})} \|Q_\veps\|^{\sfrac{4}{3}}_{L^6(\Omega)} \bigg)^{\sfrac{3}{4}} \notag \\
& \leq 
C \|Q_\veps\|_{L^\infty(0,T; H^1(\Omega))} \|\uvect_\veps\|_{L^{\sfrac{4}{3}}(0,T; W^{1,\sfrac{4}{3}}(\Omega;\R^3)) } \notag \\
& \leq C ,  
\end{align}
using the relation   $\Lambda_\veps = \frac{1}{2}(\nabla_x \uvect_\veps - \nabla^\top_x \uvect_\veps)$. 
Finally, we shall find a suitable  bound for 
\begin{align*}
\| \Gamma H[Q_\veps,c_\veps]\|_{L^{\sfrac{4}{3}}(0,T;L^{\sfrac{12}{11}}(\Omega))} .
\end{align*}
To do so, let us first recall the expression  of $H[Q,c]$ from  \eqref{formula-H}, which  is 
\begin{align*}
H[Q,c] = \Delta_x Q - \frac{1}{2} (c-c_*) Q + b \left( Q^2 - \frac{\tr (Q^2)}{3} \mathbb I_3 \right) - c_* Q \tr(Q^2) , \ \ c_* \text{ is constant}.
\end{align*}
Using the above relation, 
we compute ($\Gamma>0$ is a constant)
\begin{align}\label{Q-t-3}
&\| \Gamma H[Q_\veps,c_\veps]\|_{L^{\sfrac{4}{3}}(0,T;L^{\sfrac{12}{11}}(\Omega))} \notag \\
& \leq C \bigg(\int_0^T \|\Delta_x Q_\veps\|^{\sfrac{4}{3}}_{L^{\sfrac{12}{11}}(\Omega)}\bigg)^{\sfrac{3}{4}}   + C \|c_\veps\|_{L^\infty((0,T)\times\Omega)} \bigg(\int_0^T \|Q_\veps\|^{\sfrac{4}{3}}_{L^{\sfrac{12}{11}}(\Omega)}\bigg)^{\sfrac{3}{4}}  \notag \\
& \ \ \ + C \bigg(\int_0^T \|Q_\veps\|^{\sfrac{8}{3}}_{L^{\sfrac{24}{11}}(\Omega)}\bigg)^{\sfrac{3}{4}} +  C \bigg(\int_0^T \|Q_\veps\|^{4}_{L^{\sfrac{36}{11}}(\Omega)}\bigg)^{\sfrac{3}{4}} \notag \\
& \leq C \|\Delta_x Q_\veps\|_{L^2(0,T; L^2(\Omega))} + C \| Q_\veps\|_{L^\infty(0,T; L^2(\Omega))} + C \|Q_\veps\|^3_{L^\infty(0,T; L^4(\Omega))} \notag \\
& \leq C,
\end{align}
thanks to the estimates \eqref{bound-c-Q-eps}. 

So, gathering the estimates \eqref{Q-t-1}--\eqref{Q-t-3}, we obtain 
\begin{align}\label{bound-Q-t-eps}
\|\partial_t Q_\veps\|_{L^{\sfrac{4}{3}}(0,T; L^{\sfrac{12}{11}}(\Omega))} \leq C ,
\end{align}
where $C>0$ is independent of $\veps$. 

\vspace*{.2cm} 

In turn, by means of the uniform estimates \eqref{bound-c-Q-eps}, \eqref{bound-c-t-eps}, \eqref{bound-Q-t-eps}, and Aubin-Lions lemma, we guarantee the strong following convergence results:
\begin{align}\label{strong-conv-c-Q-eps}
c_\veps \to c \textnormal{ strongly in } L^2((0,T) \times \Omega) \ \text{ and } 
\  Q_\veps \to Q   \text{ strongly in } 
L^2(0,T; H^1(\Omega)) .  
\end{align}

\subsection{Limit in the  continuity equation} Having in mind the uniform bound \eqref{bound-grad-rho-eps} and the limits \eqref{limit-rho-eps-1}--\eqref{E15c}, we can pass to the $\veps\to 0$ limit in the approximate continuity equation \eqref{E4}, yielding  
\begin{align}\label{continuity-eq-limit-eps}
\left[ \int_\Omega  \rho \varphi  \right]_{t = 0}^{t = \tau}&= 
\int_0^\tau \int_\Omega{ \left[ \rho \partial_t \varphi + \rho \uvect \cdot \nabla_x \varphi \right] }
\dt \notag \\ 
&- \int_0^\tau \int_{\Gamma_{\out}} \varphi \rho \uvect_B \cdot \vc{n} \ {\rm d} S_x \dt - 
\int_0^\tau \int_{\Gamma_{\inn}} \varphi \rho_B \uvect_B \cdot \vc{n}  \ {\rm d}S_x \dt ,
\end{align}
for any $\vphi \in \C^1(  [0,T]\times \overline{\Omega})$, which is the weak formulation of the continuity equation \eqref{dens-eq-1} with initial condition \eqref{initial-dens} and boundary conditions \eqref{inflow-outflow}. In particular, it holds that 
\begin{align}\label{conti-distri-eps}
\partial_t \rho + \Div_x(\rho \uvect) = 0 \ \text{ in } \mathcal D^\prime((0,T)\times \Omega),
\end{align}
and subsequently, the continuity equation is also satisfied in the sense of  renormalized solution, introduced  by DiPerna and Lions \cite{Lions-Diperna}.

\subsection{Limit in the momentum equation}

For the limit passages in the approximate momentum equation \eqref{weak-momentum-delta-limit}, we first notice that the energy estimate \eqref{Energy-ineq-delta-limit}, the   superlinearity of $F^*$ (since $F$ is a proper convex function) and the de la Vallée-Poussin criterion for equi-integrability
imply the existence of $\Ss\in L^1((0,T)\times \Omega; \R^{3\times 3}_{\sym})$ such that 
\begin{align*}
\Ss_\veps \to \Ss \text{ weakly in } L^1((0,T)\times \Omega; \R^{3\times 3}_{\sym}). 
\end{align*}
Moreover, in order to pass to the limit in the pressure and the pressure potential, we need to have  strong $L^p$-convergence of the density. This can be done by manipulating the renormalized equation in Lemma \ref{EL3}. In this regard,  one may follow the same steps from the work \cite[Section 7]{Sarka-Jan-Justyna} to achieve that   
\begin{align} \label{strong-conv-rho-eps} 
\rho_\veps \to \rho \ \textnormal{ strongly in } L^p((0,T)\times \Omega) , \ \  1\leq p< \infty . 
\end{align}
In particular, due to the continuous dependence of the pressure on  density, we readily obtain 
\begin{align}
p(\rho_\veps) \to p(\rho)  \ \text{ and } \ P(\rho_\veps) \to P(\rho) \ \text{ strongly in } L^p((0,T)\times \Omega), \ \  1\leq p< \infty.
\end{align}

Next, for the limit passages to the convective term $\rho_\veps \uvect_\veps \otimes \uvect_\veps$ and the boundary terms, we need to prove the strong convergence of the velocity field. To this end, using all the uniform bounds obtained w.r.t. $\veps$ in the momentum equation yields that 
\begin{align*}
\bigg|\int_0^T \int_{\Omega} \partial_t (\rho_\veps \uvect_\veps) \cdot \boldvphi \d x \d t \bigg| 
&=
\bigg| \int_0^T \int_\Omega \left[\rho_\veps  \uvect_\veps \otimes \uvect_\veps : \nabla_x \boldvphi + p(\rho_\veps) \Div_x \boldvphi  - \Ss_\veps : \nabla_x \boldvphi \right] \d x \d t \notag \\
& \ \ \ \  - \int_0^T \int_\Omega \left[G(Q_\veps)\mathbb I_3 - \nabla_x Q_\veps \odot \nabla_x Q_\veps \right] : \nabla_x \boldvphi \d x \d t  \notag \\
& \ \ \ \  - \int_0^T \int_\Omega \left[ Q_\veps \Delta_x Q_\veps - \Delta_x Q_\veps Q_\veps + \sigma_* c^2_\veps Q_\veps   \right] : \nabla_x \boldvphi \d x \d t \notag \\ 
& \ \ \ \ - \veps \int_0^T \int_\Omega \nabla_x \rho_\veps \cdot \nabla_x \uvect_\veps \cdot \boldvphi \d x \d t \bigg|  \\
& \leq C \|\boldvphi\|_{L^\infty(0,T; X_n) } 
\end{align*}
for any test function $\boldvphi\in L^\infty(0,T; X_n)$, where $C>0$ is independent in $\veps$. Consequently, we have 
\begin{align*}
\| \partial_t (\rho_\veps \uvect_\veps)\|_{(L^\infty(0,T; X_n))^\prime} \leq C ,
\end{align*}
and equivalently,
\begin{align*}
\| \partial_t (\rho_\veps \uvect_\veps)\|_{L^1(0,T; X^\prime_n)} \leq C .
\end{align*}
One may then apply the Aubin-Lions lemma to ensure that 
\begin{align}\label{strong-conv-rho-u-eps}
\rho_\veps \uvect_\veps \to \rho \uvect  
 \ \text{ strongly in } \ L^2(0,T; X^\prime_n) . 
\end{align}
From the result \eqref{strong-conv-rho-u-eps} and the convergence of $\uvect_\veps$ in \eqref{E15b}, we can further conclude that 
\begin{align*}
\|\sqrt{\rho_\veps} \uvect_\veps\|_{L^2((0,T)\times \Omega)} \to 
\|\sqrt{\rho} \uvect\|_{L^2((0,T)\times \Omega)}.
\end{align*}
But we have, due to \eqref{bound-sqrt-rho-u-eps}, 
\begin{align*}
\sqrt{\rho_\veps} \uvect_\veps  \to 
\sqrt{\rho} \uvect \  \text{ weakly in } \ {L^2((0,T)\times \Omega)}, 
\end{align*}
which gives a  better convergence result, namely 
\begin{align*}
\sqrt{\rho_\veps} \uvect_\veps  \to 
\sqrt{\rho} \uvect \  \text{ strongly in } \ {L^2((0,T)\times \Omega)}.
\end{align*}

Then, the strong convergence of the density \eqref{strong-conv-rho-eps}, the lower bound of $\rho_\veps$  in  \eqref{bound-rho-up-low}, and the equivalence of norms in the finite dimensional space $X_n$, we infer that
\begin{align}\label{strong-conv-u-eps}
\uvect_\veps \to \uvect \ \text{ strongly in } \ L^p(0,T; X_n), 
\ \ 1\leq p< \infty . 
\end{align}

\vspace*{.2cm} 

Gathering all the convergence results, we can now pass to the $\veps \to 0$ limit in the momentum equation \eqref{weak-momentum-delta-limit}, resulting in 
\begin{align}\label{weak-momentum-eps-limit}
\int_\Omega \rho \uvect (\tau) \cdot \boldvphi(\tau) \d x &=
\int_0^\tau \int_\Omega \left[\rho \uvect \cdot \partial_t \boldvphi + \rho \uvect \otimes \uvect : \nabla_x \boldvphi + p(\rho) \Div_x \boldvphi  - \Ss : \nabla_x \boldvphi \right] \d x \dt \notag \\
& \ \ \ - \int_0^\tau \int_\Omega \left[G(Q)\mathbb I_3 - \nabla_x Q \odot \nabla_x Q \right] : \nabla_x \boldvphi \d x \d t  \notag \\
& \ \ \ - \int_0^\tau \int_\Omega \left[ Q \Delta_x Q - \Delta_x Q Q + \sigma_* c^2 Q   \right] : \nabla_x \boldvphi \d x \d t \notag \\ 
& 
\ \ \ + \int_\Omega \rho_0 \uvect_0 \cdot \vphi(0) \dx ,
\end{align} 
for any $\boldvphi \in \C^1([0,T]; X_n)$, with the initial condition
\begin{equation*}
(\rho \uvect)(0, \cdot)  =\rho_0 \uvect_0, \ \uvect_0 = \vc{v}_0 + \uvect_B, \ \vc{v}_0 \in X_n.
\end{equation*}

\subsection{Limit in the concentration and nematic order parameter equation} We recall the convergences of the quantities $c_\veps$ and $Q_\veps$ from \eqref{conv-c-Q-eps} and \eqref{strong-conv-c-Q-eps}; using these limits and the strong convergence of $\uvect_\eps$ from \eqref{strong-conv-u-eps}, we have the following limiting weak formulations for the concentration equation of active particles and the equation for $Q$ tensor parameter: 
\begin{align}\label{weak-concen-eps-limit}
\int_\Omega c(\tau)\psi (\tau) \dx - \int_0^\tau \int_\Omega c \partial_t \psi \dx\dt + \int_0^\tau \int_\Omega (\uvect \cdot \nabla_x) c \psi \dx\dt \notag 
\\ = \int_\Omega c_0 \psi(0) \dx - D_0 \int_0^\tau \int_\Omega \nabla_x c \cdot \nabla_x \psi \dx\dt 
\end{align}
for any $\tau \in [0,T]$ and any test function $\psi \in \C^1([0,T]\times \overline \Omega)$, with $c(0,\cdot)=c_0$ in $\Omega$, and 
\begin{align}\label{weak-nematic-eps-limit}
&\int_\Omega Q(\tau):\Psi (\tau) \dx - \int_0^\tau \int_\Omega Q : \partial_t \Psi \dx\dt + \int_0^\tau \int_\Omega (\uvect \cdot \nabla_x) Q : \Psi \dx\dt  \notag \\ 
& \quad + \int_0^\tau \int_\Omega [Q\Lambda - \Lambda Q]: \Psi \dx\dt
 = \int_\Omega Q_0 : \Psi(0) \dx - \Gamma \int_0^\tau \int_\Omega H[Q,c] :  \Psi \dx\dt 
\end{align}
for any $\tau \in [0,T]$ and any test function $\Psi \in \C^1_c([0,T]\times \Omega; \R^{3\times 3})$, with $Q(0,\cdot)=Q_0$ in $\Omega$. Moreover, it holds that $Q\in S^3_0$.

\subsection{Limiting energy inequality} Finally,  passing to the $\veps \to 0$  limit in the energy inequality \eqref{Main-Energy-limit-delta}, the limiting energy inequality reduces to 
\begin{align}\label{Main-energy-eps-limit}
&\frac{\d}{\d t} \int_{\Omega} \left(\frac{1}{2}\rho |\uvect - \uvect_B|^2 + P(\rho) + \frac{1}{2} |c|^2 + \frac{1}{2}|Q|^2 + \frac{1}{2}|\nabla_x Q|^2 + \frac{c_*}{4} |Q|^4 \right)   \d x  \notag \\
& + \frac{1}{2}\int_{\Omega} \left[F(\D_x\uvect) + F^*(\Ss) \right] \d x  
+  \frac{D_0}{2} \int_\Omega |\nabla_x c|^2 \d x \notag + 
\frac{\Gamma}{4} \int_{\Omega} |\Delta_x Q|^2  \d x  + 
\frac{c_*^2 \Gamma}{2} \int_{\Omega} |Q|^6\dx    \notag \\
& + \int_{\Gamma_{\out}} P(\rho) \uvect_B \cdot \mathbf n  \d S_x  + \int_{\Gamma_{\inn}} P(\rho_B) \uvect_B \cdot \mathbf n \d S_x  
\notag \\
  \leq  & C \int_{\Omega} \left( \frac{1}{2} |c|^2 + \frac{1}{2}|Q|^2 + \frac{1}{2}|\nabla_x Q|^2 + \frac{c_*}{4} |Q|^4 \right)   \d x   + \int_\Omega \Ss : \nabla_x \uvect_B \dx   \notag \\ 
  & - \int_{\Omega} \left[ \rho\uvect\otimes \uvect + p(\rho) \mathbb{I}_3\right]  : \nabla_x \uvect_B \d x  + \int_\Omega \rho \uvect \cdot (\uvect_B \cdot \nabla_x \uvect_B) \d x 
 \notag \\
 & - \frac{1}{2} \int_{\partial \Omega} \left(\frac{1}{2} |Q_B|^2 +  \frac{c_*}{4} \tr^2(Q_B^2) \right) \uvect_B \cdot \mathbf n   \d S_x 
+ 2c_* \Gamma \int_{\partial \Omega} |Q_B|^2 Q_B : \nabla_x Q_B \cdot \mathbf n \d S_x + C ,
\end{align}
where the constant $C>0$ is uniform w.r.t. $n$.

In \eqref{Main-energy-eps-limit}, one may apply the Gr\"{o}nwall's lemma to further deduce that
\begin{align}\label{Energy-ineq-eps-limit}
& \left[\int_{\Omega} \left(\frac{1}{2}\rho |\uvect - \uvect_B|^2 + P(\rho) + \frac{1}{2} |c|^2  + \frac{1}{2}|Q|^2  + \frac{1}{2}|\nabla_x Q|^2 + \frac{c_*}{4}|Q|^4 \right)   \d x\right]_{0}^{\tau}  \notag \\
& 
+  \int_0^\tau \left( \frac{D_0}{2}\|\nabla 
_x c\|^2_{L^2(\Omega)}  + \frac{\Gamma}{4}\|\Delta_x Q\|^2_{L^2(\Omega)} + \frac{c^2_*\Gamma}{2}\|Q\|^{6}_{L^6(\Omega)} \right) \d t   \notag \\
& + \frac{1}{4}\int_0^\tau \int_{\Omega}  \left[F (\D_x \uvect) + F^*(\Ss) \right] \d x \d t    + \int_0^\tau \int_{\Gamma_{\out}} P(\rho) \uvect_B \cdot \mathbf n  \d S_x \d t 
\notag \\
 &\leq  C 
\end{align}
for a.a. $\tau \in [0,T]$, where the constant $C>0$ is independent of  $n$.

\vspace*{.2cm}

Combining the above, we state the following proposition.
\begin{prop}\label{Prop-Approx-III}
Assume that the hypothesis in Proposition \ref{Proposition-Approx-I} hold true and  let $(\rho_\veps, \uvect_\veps, c_\veps, Q_\veps,\Ss_\veps)_{\veps>0}$ be the sequence of functions obtained in Section \ref{Section-delta-to-0} for given  $n\in \mathbb N$. Then, there exist some functions $(\rho, \uvect, c, Q, \Ss)$ such that the following convergence results hold:
\begin{align*}
& \rho_\veps \to \rho \, \text{ strongly in } \, L^p((0,T)\times \Omega) , \ 1 \leq p < \infty, \ \  \rho_\veps \uvect_\veps \to \rho \uvect \, \text{ strongly in } \, L^2(0,T; X^\prime_n)  , \\
& \veps \nabla_x \rho_\veps \to 0 \, \text{ weakly in } \, L^2((0,T)\times \Omega), \ \ \rho_\veps \to \rho \, \text{ weakly-}(*) \text{ in } \, L^\infty((0,T)\times \partial \Omega; \partial S_x) , \\ 
& \uvect_\veps \to \uvect \, \text{ strongly in } \, L^p((0,T; X_n), \ 1\leq p< \infty , \ \ \Ss_\veps \to \Ss \, \text{ weakly in } \, L^1((0,T)\times \Omega; \mathbb R^{3\times 3}_{\sym}) , \\   
&  c_\veps \to c \, \text{ weakly in } \, L^2(0,T; H^1(\Omega))  , \ \ c_\veps \to c \, \text{ strongly in } \, L^2((0,T)\times \Omega) , \\
&  Q_\veps \to Q \, \text{ weakly in } \, L^2(0,T; H^2(\Omega)) , \ \ Q_\veps \to Q \, \text{ strongly in } \, L^2(0,T; H^1(\Omega))  . 
\end{align*}
Moreover, $(\rho, \uvect, c, Q, \Ss)$ satisfy the weak formulations presented in this section.  
\end{prop}

\vspace*{.2cm}

 Here, we conclude with the $n$-th level solutions. The sequence of limiting functions  given by Proposition \ref{Prop-Approx-III} can be denoted  as  $(\rho_n, \uvect_n, c_n, Q_n, \Ss_n)_{n\in \mathbb N}$. We shall treat them in the $n\to\infty$ limit, in the next section.

\section{Final level: the $n\to \infty$ limit }\label{Section-limit-n}

In the final step, we pass to the $n \to \infty$ limit, that is, in the Galerkin approximation. To this end, we consider regular (as required in the Galerkin level) sequences  $\rho_{0,n}$, $(\rho \uvect)_{0,n}$, $c_{0,n}$,  $Q_{0,n}$ and $Q_{B,n}$ such that 
\begin{align*}
\rho_{0,n} \to \rho_0  \text{ in } L^\gamma(\Omega),   \ \ \ \frac{|(\rho \uvect)_{0,n}|^2}{\rho_{0,n}} \to \frac{|(\rho \uvect)_0|^2}{\rho_0} \text{ in } L^1(\Omega) , \\
c_{0,n} \to c_0 \text{ in } L^\infty(\Omega), 
\ \ \  Q_{0,n} \to Q_0 \text{ in } \text{ in } H^1(\Omega)   ,\\
 Q_{B,n} \to Q_{B} \text{ in } \text{ in } H^{\sfrac{3}{2}}(\partial \Omega).   
\end{align*}
Then, due to the  energy estimate \eqref{Energy-ineq-eps-limit}, we have 
\begin{align}\label{uniform-bound-n}
\|\sqrt{\rho_n} \uvect_n\|_{L^\infty(0,T; L^2(\Omega))} + \|P(\rho_n)\|_{L^\infty(0,T;L^1(\Omega))} + 
\|\uvect_n\|_{L^{\sfrac{4}{3}}(0,T ;  W^{1,\sfrac{4}{3}}(\Omega;\R^3)  ) } \leq C ,
\end{align}
uniformly w.r.t.  $n\to \infty$. Using the above estimates and the continuity equation \eqref{conti-distri-eps} at the $n$-th level, one can deduce, up to a subsequence, if necessary, that 
\begin{align}\label{Limit-several-n}
\begin{dcases} 
\rho_n \to \rho \ \text{ in } \  \C_{\textnormal{weak}}([0,T]; L^\gamma(\Omega)), \\ \rho_n |_{\Gamma_{\out}} \to \rho \ \text{ weakly-$(*)$ in } \ 
L^\infty(0,T; L^\gamma(\Gamma_{\out}; |\uvect_B\cdot \mathbf n| \d S_x)   , \\
\uvect_n \to \uvect \ \text{ weakly in } \ 
L^{\sfrac{4}{3}}(0,T ;  W^{1,\sfrac{4}{3}}(\Omega;\R^3)  ) ,\\
 \rho_n \uvect_n \to \mathbf m  \ \text{ weakly-$(*)$ in } \ L^\infty(0,T; L^{\frac{2\gamma}{\gamma+1}}(\Omega; \mathbb R^3))  .
\end{dcases}
\end{align}
We have to show that 
\begin{align}\label{a.e.-m-rho-u}
\mathbf m = \rho \uvect \  \text{ a.e. in } (0,T)\times \Omega . 
\end{align}
To this end, we write the following 
variant of the Div-Curl lemma,
whose proof can be found in \cite[Lemma 8.1]{Dissipative-Feireisl}. 
\begin{lemma}\label{EWANL1}
 Let $T>0$ and $\Omega \subset \mathbb R^3$ be a bounded domain.  Suppose that 
\begin{align*} 
r_n \to  r  \ \mbox{weakly in} \ L^p((0,T)\times \Omega), \ \  \mathbf v_n \to \mathbf v \ \mbox{weakly in}\ L^q(0,T; L^q( \Omega; \R^3)), \ \ p > 1, \ q > 1, 
\end{align*} 
and
\begin{align*}
r_n  \mathbf v_n \to \mathbf w \ \mbox{weakly in}\ L^r(0,T; L^r( \Omega; \R^3)), \  r >  1.
\end{align*} 
In addition, let 
\begin{align*}
\partial_t r_n + \Div_x (r_n \mathbf v_n)=0  \ \mbox{ in} \ \mathcal{D}'((0,T)\times \Omega), 
\end{align*}
and
\begin{align*}
\left\| \nabla_x \mathbf v_n \right\|_{L^1(0,T; L^1( \Omega; \R^3))} \leq C ,
\end{align*} 
for some constant $C>0$ independent in $n$. 
Then, it holds that 
\begin{align*} 
\mathbf w = r  \mathbf v  \ \mbox{ a.e. in} \  (0,T)\times \Omega . 
\end{align*}
\end{lemma}

\vspace{.2cm}

A direct application of Lemma \ref{EWANL1} to 
\begin{align*}
r_n = \rho_n,  \ \mathbf v_n = \uvect_n, \ \text{ and } \ p = \gamma, \ q = \frac{4}{3}, \  r = \frac{2\gamma}{\gamma+1} 
\end{align*}
yields \eqref{a.e.-m-rho-u}.

\subsection{Limit in the continuity equation}
The convergence results \eqref{Limit-several-n} and \eqref{a.e.-m-rho-u} allow us to pass to the limit in the continuity equation \eqref{continuity-eq-limit-eps} w.r.t. $n\to \infty$. Accordingly it holds that
\begin{align*}
\partial_t \rho + \Div_x (\rho \uvect) = 0 \ \text{ in } \mathcal D^\prime( (0,T)\times \Omega) . 
\end{align*}

\subsection{Limit in the momentum equation}\label{Section-subsub-Momentum}

For the limit passage in the momentum balance \eqref{weak-momentum-eps-limit}, we  infer from the energy estimate \eqref{Energy-ineq-eps-limit} and the uniform  bounds \eqref{uniform-bound-n} that there exist 
\begin{align*}
\mathbb{S} \in L^1((0,T)\times \Omega; \mathbb{R}^{3 \times 3}_{\text{sym}}), \quad \ 
\overline{ p(\rho)}, \ \overline{P(\rho)}, \ \overline{\rho |\uvect|^2} \in 
L^\infty (0,T;\mathcal{M}( \overline{\Omega})), 
\end{align*}
\begin{align*} 
 \overline{\rho \uvect \otimes \uvect} \in 
L^\infty (0,T; \mathcal{M} ( \overline{\Omega};\mathbb{R}^{3 \times 3}) ) , 
\end{align*}
such that
\begin{align}\label{S-P-CONVECTIVE-conv-n}
\begin{dcases} 
 \mathbb{S}_n \to \mathbb{S} \ \text{ weakly in } \
 L^1((0,T)\times \Omega; \mathbb{R}^{3 \times 3}_{\text{sym}}) ,   \\ 
 p(\rho_n) \to \overline{p(\rho)} \  
\text{ weakly-$(*)$ in } \ L^\infty(0,T;\mathcal{M}( \overline{\Omega})), \\ 
P(\rho_n ) \to \overline{P(\rho)} \ \text{ weakly-$(*)$ in } \ L^\infty(0,T;\mathcal{M}( \overline{\Omega} ) ), \\
 \rho_n |\uvect_n|^2 \to 
 \overline{\rho|\uvect|^2} \ \text{ weakly-$(*)$ in } \  L^\infty (0,T;\mathcal{M}(\overline{\Omega})),  \\
\rho_n \uvect_n \otimes \uvect_n \to  \overline{\rho \uvect \otimes \uvect} \  \text{ weakly-$(*)$ in } \ L^\infty(0,T;\mathcal{M}(\overline{\Omega}; \R^{3\times 3} )). 
\end{dcases}
\end{align}
Observe that we are not yet able to identify the limit functions in the convergences above \eqref{S-P-CONVECTIVE-conv-n}. To this end,  we introduce the error terms $\Re$, $\mathfrak{E}$, 
\begin{align}
&\Re := \Big[ \overline{\rho \uvect \otimes \uvect} + \overline{p(\rho)} \mathbb I_3 \Big] - \Big[ \rho \uvect \otimes \uvect - p(\rho) \mathbb I_3 \Big] \in L^\infty (0,T;\mathcal{M}^+( \overline{\Omega};\mathbb{R}^{3 \times 3}_{\text{sym}} )), \label{errorterm-1} 
\\
&\mathfrak{E} := \left[\frac{1}{2} \overline{\rho |\uvect|^2} + \overline{P(\rho)} \right] - \left[\frac{1}{2} \rho|\uvect|^2 + P(\rho) \right]  \in L^\infty (0,T;\mathcal{M}^+( \overline{\Omega} ) ) ,  \label{errorterm-2}
\end{align}
 where $\rho$ and $\uvect$ are the limits obtained in \eqref{Limit-several-n}.
From the convexity of the functions $P - \underline{a}p$ and $\overline{a}p - P$ assumed in \eqref{S2}, we know that
\begin{align}
\underline{a} \left( \overline{p(\rho)} - p(\rho) \right) \leq \overline{P(\rho)} - P(\rho)  \ \text{ and } \ \overline{P(\rho)} - P(\rho) \leq \overline{a} \left( \overline{p(\rho)} - p(\rho) \right). \nonumber
\end{align}
For the error terms defined in \eqref{errorterm-1}, \eqref{errorterm-2}, this implies the existence of constants $0 < \underline{d} \leq \overline{d}< \infty$ such that
\begin{align}
\underline{d} \mathfrak{E} \leq \operatorname{tr} \left[\Re \right] \leq \overline{d} \mathfrak{E} . \label{compconlim}
\end{align}

Finally, performing similar steps as for the strong convergences \eqref{strong-conv-c-Q-eps}  of the concentration of active particles and $Q$-tensor parameter, one can show that
\begin{align}\label{strong-conv-c-Q}
c_n \to c \, \text{ strongly in } L^2((0,T)\times \Omega) \  \text{ and } \  Q_n \to Q \, \text{ strongly in } L^2(0,T; H^1(\Omega)) ,  
\end{align} in the $n\to\infty$ limit.

Gathering the above convergence results and by introducing the error terms given by \eqref{errorterm-1}--\eqref{errorterm-2},  we can now pass to the limit in the momentum balance \eqref{weak-momentum-eps-limit} w.r.t. $n$.  This yields
\begin{align}\label{weak-momentum-n-limit}
\int_\Omega \rho \uvect (\tau) \cdot \boldvphi(\tau) \d x &=
\int_0^\tau \int_\Omega \left[\rho \uvect \cdot \partial_t \boldvphi + \rho \uvect \otimes \uvect : \nabla_x \boldvphi + p(\rho) \Div_x \boldvphi  - \Ss : \nabla_x \boldvphi \right] \d x \dt \notag \\
& \ \ \ - \int_0^\tau \int_\Omega \left[G(Q)\mathbb I_3 - \nabla_x Q \odot \nabla_x Q \right] : \nabla_x \boldvphi \d x \d t  \notag \\
& \ \ \ - \int_0^\tau \int_\Omega \left[ Q \Delta_x Q - \Delta_x Q Q + \sigma_* c^2 Q   \right] : \nabla_x \boldvphi \d x \d t \notag \\ 
& 
\ \ \   - \int_0^\tau \int_{\Omega} \nabla_x \boldvphi : \d \Re(t)  + \int_\Omega (\rho \uvect)_0\cdot \vphi(0)  \dx ,
\end{align} 
for any test function $\boldvphi \in \C^1([0,T]; X_N)$, for any   
$N\in \mathbb N$, chosen arbitrarily. 

It remains to show that the weak formulation \eqref{weak-momentum-n-limit} holds for general test functions $\boldvphi\in \C^1_c([0,T]\times \Omega)$ . In fact, by the choices of $X_n$ introduced in \eqref{Space-X_n}, for any $\boldvphi\in \C^1_c([0,T]\times \Omega)$, one may find a sequence $\{\boldvphi_N\}_{N\in \mathbb N}$ such that 
\begin{align*}
\boldvphi_N \to \boldvphi \, \text{ in } \, W^{1, \infty}((0,T)\times \Omega) \, \text{ as }\, N\to \infty. 
\end{align*}
Thus, we can extend the validity of the momentum balance \eqref{weak-momentum-n-limit} to any test function $\boldvphi \in \C^1_c([0,T]\times \Omega)$.

\subsection{Limit in the concentration and tensor parameter equations} 

Collecting all the required convergence results at this level,  we can pass to the limit in the concentration equation for active particles  \eqref{weak-concen-eps-limit} and the equation for $Q$-tensor parameter \eqref{weak-nematic-eps-limit} w.r.t. $n\to \infty$. The limiting formulations have the same forms \eqref{weak-concen-eps-limit} and   \eqref{weak-nematic-eps-limit}  as in the previous level;  we skip the details for brevity.

\subsection{Limiting energy inequality}
The energy inequality reduces to the following \eqref{Main-energy-eps-limit}, in the $n\to \infty$ limit:
\begin{align}\label{Main-energy-n-limit}
&\frac{\d}{\d t} \int_{\Omega} \left(
\frac{1}{2}\overline{\rho |\uvect|^2} - \rho \uvect\cdot \uvect_B + \frac{1}{2}\rho |\uvect_B|^2 + \overline{P(\rho)} + \frac{1}{2} |c|^2 + \frac{1}{2}|Q|^2 + \frac{1}{2}|\nabla_x Q|^2 + \frac{c_*}{4} |Q|^4 \right)   \d x  \notag \\
& + \frac{1}{2}\int_{\Omega} \left[F(\D_x\uvect) + F^*(\Ss) \right] \d x  
+  \frac{D_0}{2} \int_\Omega |\nabla_x c|^2 \d x \notag + 
\frac{\Gamma}{4} \int_{\Omega} |\Delta_x Q|^2  \d x  + 
\frac{c_*^2 \Gamma}{2} \int_{\Omega} |Q|^6\dx    \notag \\
& + \int_{\Gamma_{\out}} P(\rho) \uvect_B \cdot \mathbf n  \d S_x  + \int_{\Gamma_{\inn}} P(\rho_B) \uvect_B \cdot \mathbf n \d S_x  
\notag \\
  \leq  & C \int_{\Omega} \left(\frac{1}{2} |c|^2 + \frac{1}{2}|Q|^2 + \frac{1}{2}|\nabla_x Q|^2 + \frac{c_*}{4} |Q|^4 \right)   \d x   + \int_\Omega \Ss : \nabla_x \uvect_B \dx \notag \\ 
  & - \int_{\Omega} \left[ \overline{\rho\uvect\otimes \uvect} + \overline{p(\rho)} \mathbb{I}_3\right]  : \nabla_x \uvect_B \d x  + \int_{\Omega} \rho \uvect \cdot (\uvect_B \cdot \nabla_x \uvect_B) \d x 
  \notag \\
 & - \frac{1}{2} \int_{\partial \Omega} \left(\frac{1}{2} |Q_B|^2 +  \frac{c_*}{4} \tr^2(Q_B^2) \right) \uvect_B \cdot \mathbf n   \d S_x 
 + 2c_* \Gamma \int_{\partial \Omega} |Q_B|^2 Q_B : \nabla_x Q_B \cdot \mathbf n \d S_x + C ,
\end{align}
where the constant $C>0$ is independent of $\delta, \varepsilon$ and $n$.

By means  of the error terms $\Re$ and $\mathfrak{E}$ defined in \eqref{errorterm-1}--\eqref{errorterm-2}, we  can now   rewrite \eqref{Main-energy-n-limit} as follows:
\begin{align}\label{Main-energy-final}
&\frac{\d}{\d t} \int_{\Omega} \left(
\frac{1}{2}\rho |\uvect-\uvect_B|^2  + P(\rho) + \mathfrak{E} + \frac{1}{2} |c|^2 + \frac{1}{2}|Q|^2 + \frac{1}{2}|\nabla_x Q|^2 + \frac{c_*}{4} |Q|^4 \right)   \d x  \notag \\
& + \frac{1}{2}\int_{\Omega} \left[F(\D_x\uvect) + F^*(\Ss) \right] \d x  
+  \frac{D_0}{2} \int_\Omega |\nabla_x c|^2 \d x \notag + 
\frac{\Gamma}{4} \int_{\Omega} |\Delta_x Q|^2  \d x  + 
\frac{c_*^2 \Gamma}{2} \int_{\Omega} |Q|^6\dx    \notag \\
& + \int_{\Gamma_{\out}} P(\rho) \uvect_B \cdot \mathbf n  \d S_x  + \int_{\Gamma_{\inn}} P(\rho_B) \uvect_B \cdot \mathbf n \d S_x  
\notag \\
  \leq  & C \int_{\Omega} \left(\frac{1}{2} |c|^2 + \frac{1}{2}|Q|^2 + \frac{1}{2}|\nabla_x Q|^2 + \frac{c_*}{4} |Q|^4 \right)   \d x   + \int_\Omega (\Ss - \Re) : \nabla_x \uvect_B \dx \notag \\ 
  & - \int_{\Omega} \left[ \rho\uvect\otimes \uvect + p(\rho) \mathbb{I}_3\right]  : \nabla_x \uvect_B \d x  
  + \int_{\Omega} \rho \uvect \cdot (\uvect_B \cdot \nabla_x \uvect_B) \d x 
  \notag \\
 & - \frac{1}{2} \int_{\partial \Omega} \left(\frac{1}{2} |Q_B|^2 +  \frac{c_*}{4} \tr^2(Q_B^2) \right) \uvect_B \cdot \mathbf n   \d S_x 
 + 2c_* \Gamma \int_{\partial \Omega} |Q_B|^2 Q_B : \nabla_x Q_B \cdot \mathbf n \d S_x + C ,
\end{align}
where the effective viscous stress for the concerned flow can defined as 
\begin{align*}
\Ss_{\textnormal{eff}} : = \Ss - \Re ,
\end{align*}
and $\Re$ being the Reynolds stress.

\vspace{.3cm}

This basically concludes the proof of the {\em global existence of dissipative solution} (according to  Definition \ref{Definition-sol}), that is, the proof of Theorem \ref{Theorem-Main}.

\section{Further remarks}
\label{Section-Conclusion}   In this paper, we prove the existence of dissipative solutions for compressible active nematodynamics, where the viscous stress tensor satisfies an implicit rheological law \eqref{stress-tensor} w.r.t. some convex potential $F$. This is a significant advance in the theoretical underpinnings of active nematodynamics, since the passage from dissipative solutions to weak solutions and then to strong solutions needs to be understood for assessing the validity and robustness of these phenomenological non-equilibrium models. 
Note that, in \cite{Dissipative-Feireisl}, the authors work with the potential $F$ which satisfies 
\begin{align}\label{Coercivity-F-2} 
F(\D) \geq \mu_0 \left| \D - \frac{1}{3} \tr [\D] \mathbb I_3 \right|^{q} \quad \text{for all } \, \D \in \mathbb R^{3\times 3}_{\sym}  , \ \ q>1 , 
\end{align}
for some $\mu_0>0$.  
In the present work, we restrict ourselves to $q=\frac{4}{3}$,  to overcome the difficulty caused by the concentration of active particles. Indeed, there is a term  $\sigma_* \Div_x(c^2Q)$ in the momentum equation \eqref{momen-eq-1} stemming from the active stress, and we need the restriction of $q=\frac{4}{3}$ (see 
\eqref{esti-7}) to estimate this term in the energy inequality. However, in the case of passive nematodynamics with $\sigma_*=0$, one can work with more general convex potentials $F$ (in other words, viscous stress tensor $\Ss$) satisfying \eqref{Coercivity-F-2} for any $q>1$. The methods also differ from the work on weak solutions for compressible active nematics in three dimensions in \cite{C-Majumdar-W-Z-Comp}, originating from the different choices of the boundary conditions and the final $n\to\infty$ limit depends on entirely different concepts of the corrective Reynolds stress and energy terms in the framework of dissipative solutions.

Dissipative solutions are important from a numerical standpoint. In particular, this class of solutions is large enough to accommodate the consistency and stability of  numerical schemes such as  finite volume method, mixed finite volume - finite element method and finite difference method.  For example, we refer to \cite[Chapter 7.6]{Feireisl-Bangwei} for convergence studies of approximated (dissipative) solutions in the compressible Navier-Stokes framework. For more on numerical schemes for compressible fluids, please see \cite[Chapters 8 -- 14]{Feireisl-Bangwei}. 

There are several generalisations of our work on dissipative solutions for active nematodynamics, i.e., to general pressure potentials, generic rheological laws to account for large classes of complex fluids and generic boundary conditions such as mixed boundary conditions for ${Q}$ and physically relevant boundary conditions for $\uvect$, some of which may be outside the scope of traditional weak formulations. Dissipative solutions can provide a robust analytic and numerical framework for studying the well-posedness of these problems and the accuracy and convergence of associated numerical schemes, ultimately feeding into new-age applications of active nematics such as regenerative medicine, tissue engineering, drug delivery etc. There are numerous technical challenges stemming from the new physics of each problem, but our foundational work can go some way in popularising the concept of dissipative solutions for nematodynamics and setting the scene for larger-scale applications-oriented studies of active nematics in the future.

\appendix 

\section{Some auxiliary results}\label{Section-Appendix}

In this section, we provide some auxiliary results which has been intensively used in this paper to find proper a-priori estimate in Section \ref{Section-Apriori}. 
First, we report the following result from \cite[Lemma A.1]{C-Majumdar-W-Z-InComp}.
\begin{lemma}\label{lema-aux-1}
Let $Q$ and $Q^\prime$ be two $3\times 3$ symmetric matrices and let $\Lambda'=\frac{1}{2}(\nabla_x U - \nabla_x^\top U)$ for some $U\in W^{1,q}_0(\Omega; \R^3)$ for some $q>1$.  Then, 
\begin{align}
\int_{\Omega} (\Lambda' Q^\prime - Q^\prime \Lambda' ) \Delta_x Q \d x  = \int_\Omega \Div_x (Q^\prime \Delta_x Q - \Delta_x Q Q^\prime ) \cdot U \d x . 
\end{align}
\end{lemma}

\vspace*{.2cm}

Next, we write the well-known trace inequality which can be found for instance in \cite[Theorem 1.6.6]{Brenner-Scott}.
\begin{lemma}[Trace theorem]\label{Trace-them}
Suppose that the domain  $\Omega$ has Lipschitz boundary and $1\leq p\leq \infty$. Then, there exists a constant $C>0$ (may depend on $\Omega$) such that 
\begin{align*}
\|v \|_{L^p(\partial \Omega)} \leq C \|v\|^{1-1/p}_{L^p(\Omega)} \|v\|^{1/p}_{W^{1,p}(\Omega)} \quad \forall v \in W^{1,p}(\Omega). 
\end{align*}
\end{lemma}

\vspace*{.2cm} 
Let us  define the space 
\begin{align*}
E^{q}_0(\Omega; \R^{d}): = \left\{  \mathbf v \in L^q(\Omega; \mathbb R^d) \, \mid  \,  \D_x \mathbf v - \frac{1}{d}{\tr[\D_x \mathbf v]}\mathbb{I}_3 \in L^q(\Omega; \mathbb R^{d\times d}) , \ \mathbf v|_{\partial \Omega}=0 \right\} , \ \ 1< q<\infty .
\end{align*}
Note that, $E^q_0(\Omega; \mathbb R^d)$ is a Banach space associated with the norm
\begin{align*}
\|\mathbf v\|_{E^q_0(\Omega;\R^d)} := \|\mathbf v\|_{L^q(\Omega;\R^d)}   + \left\|\D_x \mathbf v - \frac{1}{d} {\tr[\D_x \mathbf v]} \mathbb{I}_3  \right\|_{L^q(\Omega; \mathbb R^{d\times d})} . 
\end{align*}
Now, we state the following version of trace-free Korn's inequality (see \cite[Eq. (1.2)]{Trace-free-Breit}, the proof of which can be found  in \cite{Trace-free-JU-G}, and  for a more general result see  \cite[Theorem 3.1]{Trace-free-Breit}):
\begin{lemma}[Trace-free Korn's inequality]\label{Lemma-Korn-trace-free}
Let $d\geq 3$ and $\Omega\subset \mathbb R^d$  be an open bounded domain. Then, there exists a constant $\kappa>0$ such that we have the following inequality
\begin{align}\label{Trace-free-Korn-ineq}
\|\nabla_x \mathbf v \|_{L^q(\Omega; \mathbb R^{d\times d})} \leq \kappa  \left\|\D_x \mathbf v - \frac{1}{d} {\tr[\D_x \mathbf v]} \mathbb{I}_3  \right\|_{L^q(\Omega; \mathbb R^{d\times d})} 
\end{align} 
for any $\mathbf v \in E^q_0(\Omega; \mathbb R^d)$. 
\end{lemma}

\vspace*{.2cm}

Finally, we recall the De La Vall\'{e}e--Poussin criterion for equi-integrability from \cite[Chapter 6, Lemma 6.4]{Pablo-Pedregal}. 

\begin{lemma}[De La Vall\'{e}e--Poussin criterion]\label{Lemma-La-Valle}
Let $\Omega\subset \mathbb R^d$ for  $d\geq 1$ be bounded. The sequence $\{f_j\}_{j\in \mathbb N}$ is sequentially weakly relatively compact in $L^1(\Omega)$ if and only if 
\begin{align*}
\sup_{j\in \mathbb N}  \int_\Omega \psi(|f_j|) \d x < \infty  ,
\end{align*}
for some continuous function $\psi: [0,\infty) \to \mathbb R$ with 
\begin{align*}
\lim_{s \to \infty} \frac{\psi(s)}{s} = \infty . 
\end{align*}
\end{lemma}

\section*{Acknowledgments}
 {\v{S}}.N.  and K.B. have been supported by the  Praemium Academiae of {\v{S}}. Ne{\v{c}}asov{\'{a}}. The Institute of Mathematics, CAS is supported by RVO:67985840. 
    The supports are gratefully acknowledged.
    A.M. is supported by the University of Strathclyde
New Professors Fund, the Humboldt Foundation, and Leverhulme Research Project grant 
RPG-2021-401.

\bibliographystyle{plain}
\bibliography{references}

\end{document}